\documentclass[12pt]{amsart}
\usepackage{epsf,amssymb}
\usepackage{colordvi}
\usepackage{graphicx,epsf}

\theoremstyle{plain}
\newtheorem{thm}{Theorem}

\newtheorem{prop}[thm]{Proposition}
\newtheorem{example}[thm]{Example}
\theoremstyle{definition}
\newtheorem*{rmk}{Remark}

\newcommand{\C}{{\mathbb{C}}}
\newcommand{\G}{{\mathbb{G}}}
\newcommand{\R}{{\mathbb{R}}}
\renewcommand{\P}{{\mathbb{P}}}
\newcommand{\Z}{{\mathbb{Z}}}
\newcommand{\Q}{{\mathbb{Q}}}

\newcommand{\cF}{{\mathcal{F}}}
\newcommand{\cM}{{\mathcal{M}}}
\newcommand{\cO}{{\mathcal{O}}}
\newcommand{\cQ}{{\mathcal{Q}}}
\newcommand{\cS}{{\mathcal{S}}}
\newcommand{\cT}{{\mathcal{T}}}

\newcommand{\cZ}{{\mathcal{Z}}}

\newcommand{\va}{{\mathbf a}}
\newcommand{\vb}{{\mathbf b}}
\newcommand{\vc}{{\mathbf c}}
\newcommand{\vj}{{\mathbf j}}
\newcommand{\vl}{{\mathbf l}}
\newcommand{\vm}{{\mathbf m}}
\newcommand{\vr}{{\mathbf r}}
\newcommand{\vdelta}{{\mbox{\boldmath$\delta$}}}
\newcommand{\sdelta}{{\mbox{\scriptsize\boldmath$\delta$}}}
\newcommand{\vzero}{{\mathbf 0}}

\newcommand{\be}{{\mathbf e}}
\newcommand{\bE}{{\mathbf E}}
\newcommand{\bff}{{\mathbf f}}
\newcommand{\bv}{{\mathbf v}}

\newcommand{\Hom}{\mathrm{Hom}}

\newcommand{\frakm}{\mathfrak{m}}

\newcommand{\lcm}{\mathop{\mathrm{lcm}}}
\newcommand{\Span}{\mathop{\mathrm{Span}}}

\begin{document}

\title{The Equivariant Chow rings of quot schemes}  

\author{Tom Braden}
\address{Department of Mathematics\\
        University of Massachusetts\\
        Amherst, MA \ 01003\\
        USA}
\email{braden@math.umass.edu}

\author{Linda Chen}
\address{Department of Mathematics\\
        The Ohio State University\\
        Columbus, OH \ 43210\\
        USA}
\email{lchen@math.ohio-state.edu}

\author{Frank Sottile}
\address{Department of Mathematics\\
         Texas A\&M University\\
         College Station\\
         TX \ 77843\\
         USA}
\email{sottile@math.tamu.edu}

\thanks{Braden supported in part by NSF grant DMS-0201823}
\thanks{Chen supported in part by NSF VIGRE grant DMS-9810750 and NSF grant
  DMS-0432701} 
\thanks{Sottile supported in part by NSF CAREER
  grant DMS-0134860 and the Clay Mathematical Institute}

\begin{abstract}
 We give a presentation for the (integral) torus-equivariant Chow ring of the
 quot scheme, a smooth compactification of the space of rational curves of
 degree $d$ in the Grassmannian.
 For this presentation, we refine Evain's extension of the method of
 Goresky, Kottwitz, and MacPherson
 to express the torus-equivariant Chow ring in terms of the torus-fixed points and 
 explicit relations coming from the geometry of families of torus-invariant curves.
 As part of this calculation, we give a complete description of the
 torus-invariant curves on the quot scheme and show that each family is a 
 product of projective spaces.
\end{abstract}

\maketitle

\section*{Introduction}

Let $k$ be an algebraically closed field 
and let $d,n,r$ be non-negative integers with $r<n$.
We study the quot scheme $\cQ_d:=\cQ_d(r,n)$ parametrizing
quotient sheaves  on $\P^1$
of the trivial vector bundle $\cO^n_{\P^1}$ which have rank $r$ and degree $d$. 
When $r>0$, this is a compactification of the space $\cM_d$ of
parametrized rational curves of degree $d$ on the Grassmannian $G(r,n)$
of $r$-dimensional quotients of $k^n$.  
Indeed, a morphism from $\P^1$ to $G(r,n)$ of degree $d$ is
equivalent to a quotient bundle $\cO^n_{\P^1}\rightarrow\cT$ of rank $r$
and degree $d$.  

Str{\o}mme \cite{St87} showed that $\cQ_d(r,n)$ is a smooth, projective, 
rational variety of dimension 
$r(n-r)+nd$. He described the decomposition into Bia{\l}ynicki-Birula 
cells induced by an action of a one-dimensional torus $T$ on $\cQ_d$,
thereby determining its Betti numbers.   
He also gave a presentation of its integral Chow ring (Theorem 5.3, {\it
loc. cit.}) in terms of generators and relations.  However, the set of
generators is far from minimal, and the relations are given by the
annihilator of a certain class, and are therefore non-explicit.  He also
gave a more elementary set of generators for its rational Chow ring.
We compute the corresponding equivariant 
classes in Sections \ref{classes} and \ref{S:ex}.

Later, 
Bertram used the geometry of $\cQ_d$ to determine the (small) quantum
cohomology ring of the Grassmannian~\cite{Bertram}.  
He used a recursive description of the boundary $\cQ_d\setminus\cM_d$ to
show that the $3$-point genus zero Gromov-Witten invariants of the
Grassmannian are equal to particular intersection numbers on $\cQ_d$.  By
studying certain types of intersections on quot schemes, he obtained a
quantum Schubert calculus for the Grassmannian.  However, he did not need
to compute the full cohomology ring of $\cQ_d$.

Our main result is a presentation
of the $T$-equivariant Chow ring $A_T^*(\cQ_d)$ for an
action of a torus $T$ on $\cQ_d$.   This  
determines the ordinary Chow ring $A^*(\cQ_d)$.  
Our presentation gives $A_T^*(\cQ_d)$ as an explicit
subring of a direct sum of polynomial rings.
When $k=\C$, the cycle class map induces an
isomorphism between Chow and cohomology rings, so our result
also determines the $T$-equivariant and ordinary cohomology rings of $\cQ_d$.

Our presentation arises from an analysis of the localization map.
When an algebraic torus $T$ acts on a smooth variety 
$X$ with finitely many fixed points $X^T$, 
the inclusion $i\colon X^T\hookrightarrow X$ induces the
localization map of (integral) equivariant Chow rings
 \begin{equation} \label{loc map}
   i^*\colon A^*_T(X) \to A^*_T(X^T)\ =\  
     \bigoplus_{p\in X^T} A^*_T(p)\,.
 \end{equation}
Each summand $A^*_T(p)$ is canonically isomorphic to the symmetric algebra
$S$ of the character group of $T$.
When $k=\C$, there is a similar localization for rational equivariant cohomology.

When $k=\C$, Chang and Skjelbred~\cite{ChSk74} showed that the image of the
localization map in cohomology is cut out by the images of the $T$-equivariant
cohomology of components of the one-skeleton of $X$ (the points
fixed by some codimension one subtorus of $T$).
These components are closures of families of $T$-invariant curves.
In particular, when $X$ has finitely many one-dimensional $T$-orbits (whose
closures are $T$-invariant curves), Goresky, 
Kottwitz, and MacPherson~\cite{GKM98} used this to describe
the image of the localization map for equivariant cohomology.  
Each $T$-invariant curve gives a 
relation, and these GKM relations cut out the image.
Brion~\cite{brion} showed that this remains true for rational equivariant Chow
rings of varieties over any algebraically closed field.\smallskip

The torus $T=T_{k^n}\times T_{\P^1}$ acts on $\cQ_d$,
where $T_{k^n}= (\G_m)^n$ acts diagonally on $k^n$ and $T_{\P^1}=\G_m$
acts primitively on $\P^1$. 
This action has finitely many fixed points, but there are 
infinitely many $T$-invariant curves.
The GKM relations remain valid, but are now insufficient to cut out
the image of $i^*$. 
There will be extra relations coming from 
connected components of the one-skeleton of $\cQ_d$.
Brion~\cite{brion} adapted the result of Chang and Skjelbred to Chow
groups, showing that the relations given by these families are sufficient
to cut out $A^*_T(X)$, rationally.  Determining these relations explicitly
is more difficult than for the GKM relations, however, and few
cases have been worked out in detail.      

One case for which the relations are known is when $X$ is a Hilbert scheme of
points on a toric surface, which has families of $T$-invariant curves.
Following a sugestion of Brion,  Evain~\cite{evain} used 
Edidin and Graham's~\cite{EG98b} version of the 
Atiyah-Bott-Berline-Vergne
localization formula for equivariant Chow groups to give relations in terms of
ideal-membership. 
The relations come from elements of the $T$-equivariant Chow rings of
the families of $T$-invariant curves.

We discuss this in Section \ref{extended GKM theory}, and
give a more explicit formula for Evain's relations when each
component $Y$ of $X^{T'}$ has smooth $T$-invariant subvarieties $Z$
whose classes $[Z]$ generate $A^*_T(Y)$. 
We derive necessary and sufficient linear relations over $\Q$ from
Evain's ideal-membership relations. 
Lastly, we show that if the $T$-weights of the tangent space at each
fixed point are not too dependent (see Theorem~\ref{T:GKM_algebraic_Z}),
then Evain's relations also determine the integral Chow ring.

All of these additional
hypotheses hold for $\cQ_d$.  
In fact, the components $Y$ for $\cQ_d$ are relatively simple: all are products of 
projective spaces.  
As a result, we obtain explicit descriptions of the equivariant Chow
ring of $\cQ_d$, both rationally and integrally. 
\smallskip

To describe the combinatorics of fixed points in $\cQ_d^T$, we 
use the following notations.
For an element $\va = (a_1, \dots, a_n) \in \Z^n$, we define
$|\va|=\sum_i a_i$.  
We use addition and subtraction on $\Z^n$ considered as
an abelian group, and denote the identity element by 
$\vzero = (0,\dots,0)$. If 
$\va\in (\Z_{\ge 0})^n$, we set $\va! = \prod_i a_i!$, 
where $0! = 1$.  Finally we use the partial order $\va \le \vb$ to
mean $a_i \le b_i$ for all $i$.

In Section~\ref{fixed_points}, we give an
explicit parametrization of the fixed point set
$\cQ_d^T$ by a set $\cF$ of triples $(\vdelta, \va, \vb)$, where 
\begin{itemize}
 \item $\vdelta\in \{0,1\}^n$ takes the
value $1$ exactly $n{-}r$ times, so that $|\vdelta| = n-r$, and
 \item $\va,\vb$ are elements of $(\Z_{\ge 0})^n$ 
   which satisfy $|\va|+|\vb|=d$ and for which $\delta_i = 0$ 
implies $a_i = 0$ and $b_i = 0$.

 \end{itemize}
The fixed points are maximally degenerate quotient sheaves supported
at $0$ and $\infty$; the data $\va$ and $\vb$ describe the structure of 
the stalks as modules over $\cO_{\P^1}$ and as representations of $T$.

Recall that $A^*_T(p) = S$, the symmetric algebra of the character group of $T$.
We have $S=\Z[\be_1,\dots,\be_n,\bff]$, where 
$\be_1,\dots,\be_n$, and  $\bff$ are dual to the obvious basis coming from the
decomposition $T = T_{k^n}\times T_{\P^1}$.
In particular, $\bff$ restricts to the identity character
on $T_{\P^1}$ and to the trivial character of $T_{\C^n}$.
We write $S_\Q$ for $S\otimes_\Z\Q$ and sometimes $S_\Z$ for $S$, 
when we wish to emphasize our ring of scalars.
We write $S^\cF$ for the set of tuples of polynomials 
$(f_{(\sdelta,\va,\vb)}\in S \mid (\vdelta,\va,\vb)\in \cF)$.
Then $S^\cF =  A^*_T( Q_d^T)$, under the identification of $\cF$ with
$Q_d^T$.  
We exhibit the image of the localization map as a subring of $S^\cF$.

In Section \ref{Tcurves} we describe a finite set of $T$-invariant
curves which span the tangent space at each fixed point.  In Section
\ref{Tfamilies}, we describe the families of $T$-invariant curves
and their closures.  Based on this description, our relations for the
image of $i^*_{\cQ_d}$ are the following:

\begin{enumerate}
\item[{I.}] For any pair $(\vdelta,\va, \vb)$, $(\vdelta',\va', \vb')\in \cF$ with 
$\vdelta = \vdelta'$, $\va = \va'$ and $\vb = \vb'$ except in positions $i$ and $j$, and
$\delta_i = \delta'_j = 1$ and $\delta_j = \delta'_i = 0$, we have
 \[
     f_{(\sdelta,\va,\vb)}\  \equiv\  
     f_{(\sdelta',\va',\vb')} \mod \be_j - \be_i +(a'_j - a_i)\bff.
 \]
  (Note that $(a'_j - a_i) = -(b'_j - b_i)$, since $a_i + b_i = a'_j + b'_j$.)

\item[{II.}] 
 \begin{enumerate}
  \item For any pair $(\vdelta,\va, \vb)$, $(\vdelta,\va', \vb)\in \cF$
    with $\va$, $\va'$ agreeing except in positions $i$ and $j$,
we have
   \[
      f_{(\sdelta,\va,\vb)}\ \equiv\ f_{(\sdelta,\va',\vb)} \mod \be_j - \be_i +(a'_j -
      a_i)\bff.
   \]
\item For any pair $(\vdelta,\va, \vb)$, $(\vdelta,\va, \vb')\in \cF$
      with $\vb$, $\vb'$ agreeing except in positions $i$ and $j$,
we have 
  \[
     f_{(\sdelta,\va,\vb)}\ \equiv\ 
     f_{(\sdelta,\va,\vb')} \mod \be_j - \be_i + (b_i - b'_j)\bff.
  \]
\item If we have $(\vdelta,\va, \vb)$, $(\vdelta,\va', \vb)$, 
      $(\vdelta, \va, \vb')\in \cF$ satisfying both of the previous conditions
      (with the same $i$ and $j$), and in addition $a'_j - a_i = b_i - b'_j$,
      then 
  \[
      Df_{(\sdelta,\va,\vb)} - Df_{(\sdelta,\va',\vb)}
         - Df_{(\sdelta,\va,\vb')} + Df_{(\sdelta,\va',\vb')}\ \equiv\
      0 \mod \be_j - \be_i+(a'_j - a_i)\bff\,
  \]
where  $D$ is differentiation in the direction of $\be_j^\vee$, the
dual basis vector to $\be_j$.

\item[(c)$'$] Under the hypotheses of II(c),
  \[
      f_{(\sdelta,\va,\vb)} - f_{(\sdelta,\va',\vb)} 
         - f_{(\sdelta,\va,\vb')} + f_{(\sdelta,\va',\vb')}\ \equiv\ 
      0 \mod (\be_j - \be_i+(a'_j - a_i)\bff)^2\,.
  \]
\end{enumerate}

\item[{III.}] For every $(\vdelta,\va,\vb)\in\cF$ and every 
$\vzero \le \vb' < \vb$,
  \[
     \sum_{\vb'\le \vc \le \vb}  \frac{(-1)^{|\vc|}}{(\vb-\vc)!(\vc - \vb')!} 
      D^{|\vb|-|\vb'|-1}f_{(\sdelta,\va+\vc,\vb-\vc)} 
\equiv 0 \mod \bff\,,
  \]
    where $D$ is differentiation in the direction of $\bff^\vee$, the dual 
basis vector to $\bff$.
 
\item[{III$'$.}] For every $(\vdelta,\va,\vb)\in\cF$, $\vb \ne \vzero$,
\[ 
\sum_{0\leq \vc\leq \vb} \frac{(-1)^{|\vc|}}{\vc!(\vb-\vc)!} \, 
f_{(\sdelta,\va+\vc,\vb-\vc)} \equiv 0 \mod  \bff^{|\vb|}.
\]

\end{enumerate}

Relations I,  II(a), and II(b)
are standard GKM relations, while the rest come from 
families of $T$-invariant curves.    
In particular, relations II(c)/II(c)$'$ (respectively III/III$'$) 
come from certain families whose closures are isomorphic to $\P^1\times\P^1$ 
(respectively arbitrary products of projective spaces).  

\begin{thm} \label{Main_Theorem}
The rational equivariant Chow ring $A^*_T(\cQ_d)_\Q$ is isomorphic
to the set of tuples $f = (f_{(\sdelta,\va,\vb)}) \in S^\cF_\Q$
subject to the relations {\rm I}, {\rm II(a)(b)(c)}, and\/ {\rm III}.
  
The integral equivariant Chow ring $A^*_T(\cQ_d)$ is isomorphic
to the set of tuples $f = (f_{(\sdelta,\va,\vb)}) \in S^\cF_\Z$
subject to the relations {\rm I}, {\rm II(a)(b)(c)$'$}, and\/ {\rm III}$'$.  
\end{thm}

We prove Theorem \ref{Main_Theorem} in Section \ref{main proof}.
Since the equivariant Chow ring determines the ordinary 
Chow ring for smooth spaces, this gives in principle a complete
description of the Chow ring of $\cQ_d$.
%
%
The resulting computation of Betti numbers is the same as Stroome's
 computation.

The quot scheme $\cQ_d$ represents the functor which associates to a
scheme $X$ the set of flat families over $X$ of sheaves on $\P^1$ which are
quotients of $\cO_{\P^1}^n$ of rank $r$ and degree $d$.
Thus, there is a universal exact sequence
 \begin{equation}\label{E:definition}
    0\ \rightarrow\ \cS\ \rightarrow\
   \cO^n_{\P^1\times\cQ_d}\ 
   \rightarrow\  \cT\ \rightarrow 0
 \end{equation}
of sheaves on $\P^1\times\cQ_d$ where 
$\cT$ has rank $r$ and degree $d$.
This sequence of sheaves is flat, and flatness implies that 
$\cS$ is a vector bundle of degree $-d$ and rank $n{-}r$.
In particular, points of $\cQ_d$ are exact sequences of such sheaves on $\P^1$.

Str{\o}mme's generators of the rational cohomology ring of $\cQ_d$ were
K\"{u}nneth components of the Chern classes of the tautological vector bundle 
$\cS$ on $\P^1\times\cQ_d$. 
In Section~\ref{classes}, we describe the equivariant
Chern classes of $\cS$ in $A^*_T(\cQ_d)$ and thus lifts of Str\o mme's
generators to the equivariant Chow groups. 
In Section~\ref{S:ex} we work this out explicitly for
$\cQ_2(0,2)$, using Theorem~\ref{Main_Theorem} to describe
the equivariant and ordinary Chow rings
and giving explicit lifts of Str{\o}mme's generators as localized classes. 

\section{Torus-fixed points of $\cQ_d$} \label{fixed_points}

The free module $\cO^n_{\P^1}=k^n\otimes\cO_{\P^1}$ has  
basis $\be_1,\dotsc,\be_n$.  
Write $T_{k^n}$ for the group of diagonal matrices in this basis.  
Let $[x,y]$ be
coordinates on $\P^1$ with $x$ vanishing at 0 and $y$ at $\infty$.
For $T_{\P^1}$ acting on $\P^1$ with fixed
points $0$ and $\infty$, the torus $T:= T_{k^n}\times T_{\P^1}$ acts
on $\cQ_d$ naturally as indicated by the given splitting.  

The $T$-fixed points 
are indexed by triples $(\vdelta,\va,\vb)$ in the
set $\cF$ of the Introduction.
The fixed point corresponding to $(\vdelta,\va,\vb)$ is the sequence of sheaves
on $\P^1$
 \[
   \cS_{(\sdelta,\va,\vb)}\ \longrightarrow\ \cO^n_{\P^1}\ 
   \longrightarrow\ \cT_{(\sdelta,\va,\vb)}
 \]
where $\cS_{(\sdelta,\va,\vb)}$ is the image of the map 
 \[
    \bigoplus_{i=1}^n\cO_{\P^1}(-a_i-b_i)\ 
    \xrightarrow{\ \mathrm{diag}(\delta_i x^{a_i} y^{b_i} )\ }\
    \cO^n_{\P^1}\,. 
\]
We identify this fixed point with the subsheaf
$\cS_{(\sdelta,\va,\vb)}$ of $\cO^n_{\P^1}$.

We introduce the following notation.
For natural numbers $a,b$, let $\cS_{a,b}$ be the subsheaf of $\cO_{\P^1}$
which is the image of the map
\[
\cO_{\P^1}(-a-b)\ \xrightarrow{\ x^ay^b\ }\ \cO_{\P^1}\,.
\]
Under the identification of modules over $\cO_{\P^1}$ with saturated graded
modules of the homogeneous coordinate ring $k[x,y]$,  $\cS_{a,b}$ is the ideal
of $k[x,y]$ generated by $x^ay^b$. 
The quotient $\cO_{\P^1}/\cS_{a,b}$ is the skyscraper sheaf 
 \[
  \cT_{a,b}\ =\ \cO_{\P^1}/\frakm_0^a\ \oplus
  \cO_{\P^1}/\frakm_\infty^b
 \]
on $\P^1$ supported at $0$ and at $\infty$.
Here, $\frakm_p$ is the sheaf of ideals cutting out the point $p\in\P^1$.
Then we have
 \begin{eqnarray*}
   \cS_{(\sdelta,\va,\vb)} &=& \bigoplus_{\delta_i=1}
        \cS_{a_i,b_i}\cdot\be_i\,,\qquad
        \mbox{\ and}\\
   \cT_{(\sdelta,\va,\vb)} &=& \bigoplus_{\delta_j=1} \cT_{a_j,b_j}\cdot\be_j\ \oplus\  
            \bigoplus_{\delta_j=0} \cO_{\P^1}\cdot\be_j\,.
 \end{eqnarray*}
In the sum,  $\delta_j=0$ means those $j$ in
$\{1,\dotsc,n\}$ with $\delta_j=0$, and the same for $\delta_j=1$.

The tangent space to $\cQ_d$ at this fixed point
is $\Hom(\cS_{(\sdelta,\va,\vb)},\,\cT_{(\sdelta,\va,\vb)})$.  
Let $\hat{\be}_1,\dotsc,\hat{\be}_n$ be the basis dual to
$\be_1,\dotsc,\be_n$.  
For each $i,j$, set
$\bE_{ij}:=\hat{\be}_i\otimes\be_j\in \Hom(k^n,k^n)$.  
These $\bE_{ij}$ form a basis for $\Hom(k^n,k^n)$.

\begin{thm}
\label{tangentspace}
The tangent space $T_{(\sdelta,\va,\vb)}\cQ_d$ is canonically
identified with
 \begin{equation}\label{Eq:Tangent_Space}
   \bigoplus_{\delta_i=1}\bigoplus_{\delta_j=0}
     \Hom(\cS_{a_i,b_i},\, \cO_{\P^1} )\cdot \bE_{ij} \quad \oplus\quad
   \bigoplus_{\delta_i=1}\bigoplus_{\delta_j=1}
     \Hom(\cS_{a_i,b_i},\,\cT_{a_j,b_j})\cdot \bE_{ij}\ .
 \end{equation}
\end{thm}

We now give $T$-bases (bases of $T$-eigenvectors) for these summands and
determine the corresponding weights.  
Fix a basis for the character group of $T$ as follows.
Extend the action of $T_{k^n}$ on $\cO^n_{\P^1}$ to 
$T$ by letting the factor $T_{\P^1}$ act trivially.  Then we
abuse notation and denote the character of $T$ acting on the
$i$th basis vector $\be_i$ by the same symbol $\be_i$.  Thus
the dual basis element $\hat{\be}_i$ has $T$-weight $-\be_i$.

Similarly, extend the action of $T_{\P^1}$ on $\P^1$ to
an action of 
$T$ by letting $T_{k^n}$ act trivially, and
denote by $\bff$ the character of $T$ corresponding to
the action on the dense orbit $\P^1 \setminus \{0,\infty\}$.
More precisely, we can let $T_{\P^1} \cong \C^*$ act on the
homogeneous coordinates
$k[x,y]$ of $\P^1$ by $q\cdot x = qx$ and $q\cdot y = y$.
Thus $T$ acts on the rational function $z:=x/y$ with
weight $\bff$, and on the monomial $z^a=x^ay^{-a}$ with weight $a\,\bff$.

The first sum of~\eqref{Eq:Tangent_Space} involves spaces 
of the form $\Hom(\cS_{a,b},\cO_{\P^1})=H^0(\cS^*_{a,b})$.
This space of sections has a monomial $T$-basis
 \[
   \{z^c\mid -a \le c \le b\}\,.
 \]
For example, $H^0(\cS^*_{1,2})=k\cdot\{z^{-1}, z^0=1, z, z^2\}$.
Thus if $\delta_i=1$ and $\delta_j=0$, then the piece
$\Hom(\cS_{a_i,b_i},\cO_{\P^1})\cdot\bE_{ij}$ of the tangent space
has a monomial $T$-basis
$z^c\cdot\bE_{ij}$ for all $-a_i\leq c\leq b_i$.
The basis element $z^c\cdot\bE_{ij}$ has $T$-weight 
 \begin{equation}\label{Eq:TI}
    \be_j-\be_i+ c\,\bff \,.
 \end{equation}
%
%

The second sum of~\eqref{Eq:Tangent_Space} involves spaces of the form
$\Hom(\cS_{a,b},\cT_{\alpha,\beta})$. 
  Since 
$\cT_{\alpha,\beta}$ is a skyscraper sheaf supported at $0$ and
$\infty$, a map $\phi\in\Hom(\cS_{a,b},\cT_{\alpha,\beta})$ is determined  
by its actions at $0$ and at $\infty$.
At $0$, $z = x/y$ is a local parameter, so the map $\phi$
becomes
  \[
    \phi\ \colon\ z^a\C[z] \rightarrow \C[z]/{\langle z^\alpha \rangle}\,,
  \]
and thus has the form $z^{-a}f(z)$ where $f(z)$ has degree
less than $\alpha$.
At $\infty$, $z^{-1}$ is a local parameter, and
the map $\phi$ has the form $z^{b}g(z^{-1})$ where $g(z)$ is a
polynomial of degree less than $\beta$.
Thus $\Hom(\cS_{a,b},\cT_{\alpha,\beta})$ has the monomial $T$-basis
\[
   k\cdot \{z^{(\alpha-c)-a}\mid 1\leq c\leq\alpha\}\; \oplus \;
   k\cdot \{z^{b-(\beta-c)}\mid  1\leq c\leq\beta\}\,,
\]
where elements in the first summand act by zero on the stalk at
$\infty$, and elements of the second summand act by zero at $0$.

Thus, if $\delta_i=\delta_j=1$, then the summand
$\Hom(\cS_{a_i,b_i},\cT_{a_j,b_j})\cdot\bE_{ij}$  
of \eqref{Eq:Tangent_Space} has a monomial $T$-basis
 \[
   \{ z^{(a_j-c)-a_i}\cdot\bE_{ij} \mid 1\leq c\leq a_j \}\, \cup \,
   \{ z^{(b_j-c)-b_i}\cdot\bE_{ij} \mid 1\leq c\leq b_j\}
 \]
with corresponding $T$-weights
 \begin{equation}\label{Eq:TII_III}
   \be_j-\be_i + ((a_j-c) - a_i)\bff \quad\mbox{ and }\quad 
   \be_j-\be_i + (b_i - (b_j-c))\bff \,.
 \end{equation}
%
%

We note that this discussion gives a basis for $T_{(\sdelta,\va,\vb)}\cQ_d$
consisting of
 \begin{eqnarray*}
  \sum_{\delta_j=0}\sum_{\delta_i=1} (a_i+b_i+1)
  \quad +\quad  
  \sum_{\delta_i=1}\sum_{\delta_j=1} (a_j+b_j)
  &=& (n-r)\cdot(d+r) + r\cdot d\\
  &=& r\cdot(n-r) + nd\ =\ \dim\cQ_d
 \end{eqnarray*}
elements, which shows that $\cQ_d$ is smooth at the $T$-fixed point
$\cS_{(\sdelta,\va,\vb)}$, and hence everywhere. 
In the next section, we will describe $T$-invariant curves
in $\cQ_d$ incident on $\cS_{(\sdelta,\va,\vb)}$ whose
tangent directions at $\cS_{(\sdelta,\va,\vb)}$ coincide with this given
$T$-basis of $T_{(\sdelta,\va,\vb)}\cQ_d$.


\begin{example}\label{example_one}
{\rm 
 The  quot scheme $\cQ:=\cQ_2(0,2)$ of rank 0 and degree 2
 quotients of $\cO^2_{\P^1}$ has dimension $r(n{-}r)+dn=2\cdot 0+2\cdot 2=4$.
 Note that the associated Grassmannian is a point. 
 Since $r=0$, the index $\vdelta$ is the same for each fixed point,
 $\delta_1=\delta_2=1$,  and so the fixed points are indexed by quadruples
 $(a_1,a_2,b_1,b_2)$ of non-negative integers whose sum is 2. 
 Thus there are ten fixed points. 
 We represent the fixed point $(a_1,a_2,b_1,b_2)$ by two
 columns of boxes superimposed on a horizontal line where the $i$th
 column has $a_i$ boxes above the horizontal line and $b_i$ boxes below it.
 For example, 
\[
  \raisebox{-7.5pt}{\includegraphics[width=25pt]{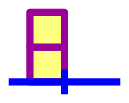}}\ 
      \leftrightarrow\ (2,0,0,0)\,,\quad
  \raisebox{-5pt}{\includegraphics[width=25pt]{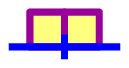}}\ 
      \leftrightarrow\ (1,1,0,0)\,,
   \qquad\mbox{and}\qquad
  \raisebox{-7pt}{\includegraphics[width=25pt]{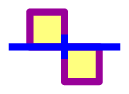}}\ 
      \leftrightarrow\ (1,0,0,1)\,.
\]

 The fixed point corresponding to 
 \raisebox{-5pt}{\includegraphics[width=20pt]{figures/11.20.00.eps}}
 is the exact sequence of sheaves on $\P^1$
\[ 
   \begin{array}{ccccccccccc}
    &&\cO(-2)\cdot\be_1&\xrightarrow{\ x^2\ }&\cO\cdot\be_1
      &\longrightarrow &\cO/\frakm_0^2\cdot\be_1\\
   0&\longrightarrow &\oplus&&\oplus&&\oplus&\longrightarrow &0\,,\\
    &&\cO\cdot\be_2&\xrightarrow{\ \,1\,\ }&\cO\cdot\be_2&\longrightarrow& 0
   \end{array}
\]
 where $\cO=\cO_{\P^1}$.
 The tangent space at this fixed point is the sum of the two 2-dimensional
 $T$-invariant spaces of homomorphisms having the indicated $T$-bases.
 \begin{eqnarray*}
    \Hom(\cO, \cO/\frakm_0^2)\cdot\bE_{21}& =&
     k\{1, z\}\cdot\bE_{21}\\
    \Hom(x^2\cO, \cO/\frakm_0^2)\cdot\bE_{11}& =&
    k\{z^{-2},z^{-1}\}\cdot\bE_{11}\
 \end{eqnarray*}
 As before, $z:=x/y$ is a local parameter at $0$ and $z^{-1}$ is a local
 parameter at $\infty$.
 These basis elements have four distinct $T$-weights
\[
   \be_1-\be_2,\ \be_1-\be_2+\bff,\ \mbox{and}\  -\bff,\ -2\bff\,.
\]

 The fixed point corresponding to 
 \raisebox{-5pt}{\includegraphics[width=20pt]{figures/11.10.01.eps}}
 is the exact sequence of sheaves on $\P^1$,
\[ 
   \begin{array}{ccccccccccc}
    &&\cO(-1)\cdot\be_1&\xrightarrow{\ x\ }&\cO\cdot\be_1
      &\longrightarrow&\cO/\frakm_0\cdot\be_1\\
   0&\longrightarrow&\oplus&&\oplus&&\oplus&\longrightarrow&0\,.\\
    &&\cO(-1)\cdot\be_2&\xrightarrow{\ y\ }&\cO\cdot\be_2
      &\longrightarrow&\cO/\frakm_\infty\cdot\be_2
   \end{array}
\]
 The tangent space at this fixed point is the sum of four 1-dimensional
 $T$-invariant spaces of homomorphisms having bases and weights as
 indicated.
\[
  \begin{array}{|l|l|r|}\hline
    \mbox{$T$-eigenspace} & \mbox{basis} &\mbox{$T$-weight}\\ \hline
    \Hom(x\cO,\cO/\frakm_0)\cdot\bE_{11}&  z^{-1}\cdot\bE_{11}&
             -\bff\ \raisebox{-3pt}{\rule{0pt}{15pt}} \\ \hline
    \Hom(x\cO,\cO/\frakm_\infty)\cdot\bE_{12}& 1\cdot\bE_{12}&
            \be_2-\be_1 \raisebox{-3pt}{\rule{0pt}{15pt}}  \\\hline
    \Hom(y\cO,\cO/\frakm_0)\cdot\bE_{21}& 1\cdot\bE_{21}&
              \be_1-\be_2 \raisebox{-3pt}{\rule{0pt}{15pt}} \\ \hline
    \Hom(y\cO,\cO/\frakm_\infty)\cdot\bE_{22}& z\cdot\bE_{22}&
              \bff\ \raisebox{-3pt}{\rule{0pt}{15pt}} \\\hline
  \end{array}
\]

}
\end{example}

\section{$T$-invariant curves} 
\label{Tcurves}

We describe a collection of $T$-invariant curves whose
tangent directions at a $T$-fixed point $(\vdelta,\va,\vb)$ form a basis for
$T_{(\sdelta,\va,\vb)}\cQ_d$.
A $T$-invariant curve on $\cQ_d$ is, by definition, a 
flat family $\cS\to\P^1$ with a $T$-action whose
fibre $\cS(s,t)$ over a point $[s,t]\in\P^1$ is a free subsheaf of
$\cO^n_{\P^1}$ of rank $n-r$ and degree $-d$.
We exhibit each curve as a subsheaf $\cS$ of $\cO^n_{\P^1\times\P^1}$ with a
$T$-action that has rank $n-r$ and degree $(-d,-1)$ with
$\cS(1,0)=\cS_{(\sdelta,\va,\vb)}$ and $\cS(0,1)=\cS_{(\sdelta',\va',\vb')}$, where
$(\vdelta',\va',\vb')$ is some other $T$-fixed point on $\cQ_d$.
Then $\cS$ defines a $T$-invariant curve on $\cQ_d$ connecting the two fixed
points. 

Each sheaf $\cS$ has one of three types: I, II, or III.
We describe them below and then argue that they have the desired
properties.  
We write $\cO$ for $\cO_{\P^1\times\P^1}$.
We use the correspondence between sheaves over $\cO$ and saturated modules over
the bihomogeneous coordinate ring $k[x,y][s,t]$ of
$\P^1_{[x,y]}\times\P^1_{[s,t]}$.
Then $\cO^n$ is the rank $n$ free module with basis $\be_1,\be_2,\dotsc,\be_n$.
\medskip

\noindent{\bf Type I:}
Let $(\vdelta',\va',\vb')$ be another fixed point where the data
$(\vdelta,\va,\vb)$ and $(\vdelta,\va',\vb')$ agree except in positions $i$ and
$j$, with $\delta_i=\delta'_j=1$ and $\delta_j=\delta'_i=0$.
Note that $a_i+b_i=a'_j+b'_j$.  
Let $\cS$ be the subsheaf of $\cO^n$ which agrees 
with both $\cS_{(\sdelta,\va,\vb)}$ and $\cS_{(\sdelta',\va',\vb')}$
except for its
component in  $\cO\cdot\be_i+\cO\cdot\be_j$, where it is the rank 1 and
degree $-(a_i{+}b_i)$ subsheaf  generated by the single element
 \begin{equation}\label{Eq:typeI}
   s x^{a_i}y^{b_i}\cdot\be_i + t x^{a'_j}y^{b'_j}\cdot\be_j \, .\smallskip
 \end{equation}

\noindent{\bf Type II:}
Let $(\vdelta',\va',\vb')$ be another $T$-fixed point where
$\vdelta=\vdelta'$, $\vb=\vb'$, and the data $\va$ and $\va'$ agree except in
positions $i$ and $j$ with $i\neq j$.
We suppose that $i$ and $j$ have been chosen so that $a_i<a'_i$.
Then $a_j>a'_j$ and $c:=a'_i-a_i=a_j-a'_j>0$.
Set $\gamma:=a_i+b_i+c-a_j-b_j$.

If $\gamma\geq 0$, let $\cS$ be the subsheaf of $\cO^n$ which 
agrees with both $\cS_{(\sdelta,\va,\vb)}$ and $\cS_{(\sdelta',\va',\vb')}$
except for its components in $\cO\cdot\be_i+\cO\cdot\be_j$,
where it is the subsheaf
generated by 
 \begin{equation}\label{Eq:typeII}
   x^{a_i+c}y^{b_i}\cdot\be_i\,,\quad
   x^{a_j}y^{b_j}\cdot\be_j\,,\quad \mbox{ and }\quad
   s x^{a_i}y^{b_i}\cdot\be_i + t x^{a_j-c}y^{b_j+\gamma}\cdot\be_j \, .
 \end{equation}
When $\gamma<0$, replace the third generator by 
$s x^{a_i}y^{b_i-\gamma}\cdot\be_i + t x^{a_j-c}y^{b_j}\cdot\be_j$.

The remaining sheaves of type II are obtained by interchanging the roles of
$\va$ and $\vb$.
That is, $\vdelta=\vdelta'$, $\va=\va'$ and $\vb,\vb'$ agree except in positions 
$i\neq j$, and we further have that $c:=b'_i-b_i=b_j-b'_j>0$.
Set $\gamma:=a_i+b_i+c-a_j-b_j$.
If $\gamma\geq 0$, let $\cS$ be as before, except with the
generators~\eqref{Eq:typeII} replaced by 
 \begin{equation}\label{Eq:typeII-b}
   x^{a_i}y^{b_i+c}\cdot\be_i\,,\quad
   x^{a_j}y^{b_j}\cdot\be_j\,,\quad\mbox{ and }\quad
   s x^{a_i}y^{b_i}\cdot\be_i + t x^{a_j+\gamma}y^{b_j-c}\cdot\be_j \, .
 \end{equation}
If $\gamma<0$, then the third generator will be 
$s x^{a_i-\gamma}y^{b_i}\cdot\be_i + t x^{a_j}y^{b_j-c}\be_j$.
\smallskip

\noindent{\bf Type III:}
Let $(\vdelta',\va',\vb')$ be another $T$-fixed point where
$\vdelta=\vdelta'$ and the data $(\va,\vb)$ and $(\va',\vb')$ agree except in
position $i$.  Thus $\delta_i\neq 0$ and $a_i+b_i=a'_i+b'_i$.
Let $\cS$ be the subsheaf of $\cO^n$ which 
agrees with both $\cS_{(\sdelta,\va,\vb)}$ and
$\cS_{(\sdelta',\va',\vb')}$, except for
its component in $\cO\cdot\be_i$, where it is the 
rank 1 and degree $-(a_i+b_i)$ subsheaf generated by
 \begin{equation}\label{Eq:typeIII}
   s x^{a_i}y^{b_i}\cdot \be_i + t x^{a'_i}y^{b'_i}\cdot\be_i\,.\smallskip
 \end{equation}

\begin{thm}\label{Th:Curves_are_correct}
 The subsheaves $\cS$ of\/ $\cO^n_{\P^1\times\P^1}$ of types {\rm I}, 
 {\rm II}, and\/ {\rm III} are $T$-invariant
 and free of rank $n{-}r$ and degree $(-d,-1)$.
 They satisfy  $\cS(1,0)=\cS_{(\sdelta,\va,\vb)}$ and 
 $\cS(0,1)=\cS_{(\sdelta',\va',\vb')}$, and hence define $T$-invariant
 curves on $\cQ_d$.
\end{thm}

\begin{proof}
  The generators of $\cS$ are $T$-invariant, except for those described
 by~\eqref{Eq:typeI}, \eqref{Eq:typeII}, \eqref{Eq:typeII-b},
 and~\eqref{Eq:typeIII}.  
 But $T$ acts transitively on those generators for $s\cdot t\neq0$.
Therefore, each sheaf $\cS$ is $T$-invariant.
In all cases, $\cS$ has degree $-1$ with respect to $\P^1_{[s,t]}$.

 The theorem is clear for the sheaves of types I and
 III, as they are constant on $\P^1_{[s,t]}$, except for the rank 1
 components~\eqref{Eq:typeI} and~\eqref{Eq:typeIII}, each of which has
 degree $(-(a_i+b_i),-1)$.
 Specializing these generators at $[s,t]=[1,0]$ and
 $[0,1]$ shows that $\cS(1,0)=\cS_{(\sdelta,\va,\vb)}$ and 
 $\cS(0,1)=\cS_{(\sdelta',\va',\vb')}$.

 We use a Gr\"obner basis argument for the sheaves of type II.
 The Hilbert function for a submodule $M$ of $\cO^2$ equals the Hilbert function
 for the module of leading terms of any Gr\"obner basis of $M$.
 As explained in~\cite[Ch.~15]{E95}, a weight $\omega$ selecting these leading
 terms induces a $\G_m$-action on $\cO^2$ whose 
 restriction to the Gr\"obner basis of $M$ generates a flat
 family over $\mathbb{A}^1$ whose special fibre is the module of
 leading terms. 

 For now, set $s=t=1$.
 If $\be_i>\be_j$, then the 
 generators~\eqref{Eq:typeII} form a Gr\"obner basis for any
 position-over-monomial ordering, and the third generator
 has leading term $x^{a_i}y^{b_i}\cdot\be_i$.
 As $c>0$, the module of leading terms is generated by
 $x^{a_i}y^{b_i}\cdot\be_i$ and $x^{a_j}y^{b_j}\cdot\be_j$, and so
 it has rank 2 and degree $-(a_i{+}b_i{+}a_j{+}b_j)$.
 The weight $\omega$ with $\omega(\be_i)=0$ and $\omega(\be_j)=-1$ induces the
 leading terms and has corresponding $\G_m$-action
 $t.(\be_i,\be_j)=(\be_i,t\be_j)$, for $t\in\G_m$.
 This action on $\cS(1,1)$ is the flat family of modules over $\mathbb{A}^1$
 generated by 
\[
   x^{a_i+c}y^{b_i}\cdot\be_i\,,\quad
   x^{a_j}y^{b_j}\cdot\be_j\,,\quad\mbox{and}\quad 
   x^{a_i}y^{b_i}\cdot\be_i + t x^{a_j-c}y^{b_j+\gamma}\cdot\be_j \ ,
\]
 which is just the part of $\cS$ in
 $\cO\cdot\be_i+\cO\cdot\be_j$ restricted to the affine subset $U$ of
 $\P^1_{[s,t]}$ where $s\neq 0$.
 Thus $\cS|_U$ is a flat family over $U$ of free subsheaves of $\cO^n_{\P^1}$
 of rank $n{-}r$ and degree $-d$, and $\cS(1,0)=\cS_{(\sdelta,\va,\vb)}$.

 When $s=t=1$ the generators~\eqref{Eq:typeII} form a Gr\"obner basis when
 $\be_i<\be_j$, where the third generator has leading term
 $x^{a_j-c}y^{b_j+\gamma}\cdot\be_j$. 
 The module of leading terms is generated by 
\[
   x^{a_i+c}y^{b_i}\cdot\be_i\,,\quad
   x^{a_j}y^{b_j}\cdot\be_j\,,\quad\mbox{and}\quad 
   x^{a_j-c}y^{b_j+\gamma}\cdot\be_j \ .
\]
 Since $\gamma\geq 0$ and $c>0$, saturating the ideal of $k[x,y]$ 
 generated by $x^{a_j}y^{b_j}$ and $x^{a_j-c}y^{b_j+\gamma}$ by the irrelevant
 maximal ideal generated by $x$ and $y$ gives
 the ideal generated by $x^{a_j-c}y^{b_j}$. 
 Thus the module of leading terms is generated by 
\[
   x^{a_i+c}y^{b_i}\cdot\be_i\,,\quad\mbox{and}\quad 
   x^{a_j-c}y^{b_j}\cdot\be_j \ .
\]
 As before, restricting $\cS$ to the affine set of points $[s,t]$ of $\P^1$ where 
 $t\neq 0$ gives a flat family of subsheaves of $\cO^n_{\P^1}$ of rank $n{-}r$
 and degree $-d$ with special fibre $\cS(0,1)=\cS_{(\sdelta',\va',\vb')}$. 
 The same arguments suffice for the module generated by~\eqref{Eq:typeII-b}.
\end{proof}

%
%
%
\begin{thm} \label{tangentdirections}
  For any $T$-fixed point  $\cS_{(\sdelta, \va, \vb)}$ in $\cQ_d$, 
  the set of tangent directions to the $T$-invariant curves induced 
  by the sheaves $\cS$ of types {\rm I}, {\rm II}, and\/ {\rm III}
  corresponds to the $T$-basis of $T_{(\sdelta, \va, \vb)}\cQ_d$ defined
  in Section~$\ref{fixed_points}$, and this correspondence is $T$-equivariant,
  respecting the weights.
  More specifically, at the $T$-fixed point $\cS_{(\sdelta, \va, \vb)}$,
\begin{enumerate}
 \item[I.]  The weight of the type {\rm I} curve~$\eqref{Eq:typeI}$ is
     \[    \be_j-\be_i+(a'_j-a_i)\bf\,,\]
     and such curves correspond to the first summand of~$\eqref{Eq:Tangent_Space}$.

 \item[II.] The weight of the type {\rm II} curve~$\eqref {Eq:typeII}$ is 
    \[   \be_j-\be_i + (a'_j-a_i)\bf\,,\]
     the weight of the type  {\rm II} curve~$\eqref {Eq:typeII-b}$ is 
    \[   \be_j-\be_i + (b_i-b'_j)\bf\,,\]
      and such curves correspond to the second summand
      of~$\eqref{Eq:Tangent_Space}$ when $i\neq j$.
  
 \item[III.]  The weight of the type {\rm III} curve~$\eqref{Eq:typeIII}$ is
     \[    (a'_j-a_i)\bf\,,\]
      and such curves correspond to the second summand
      of~$\eqref{Eq:Tangent_Space}$ when $i= j$.

\end{enumerate}
  
\end{thm}

\begin{proof}
 In each of these curves, the $T$-fixed point $\cS_{(\sdelta,\va,\vb)}$ is the point
 where $s=1$ and $t=0$.
 In what follows, we work locally, setting $s=1$.
 
 Note that the generator~\eqref{Eq:typeI} of a type I sheaf may be rewritten
\[
    x^{a_i}y^{b_i}\cdot
    ( \be_i + t\cdot x^{a'_j-a_i}y^{b'_j-b_i}\bE_{ij}(\be_i))\,.
\]
 This shows that the tangent space to the curve when $t=0$ is the $T$-basis
 element $x^{a'_j-a_i}y^{b'_j-b_i}\bE_{ij}$ of 
 $\Hom(\cS_{a_i,b_i}\cdot\be_i,\, \cO_{\P^1}\cdot\be_j)$.
 Thus the tangent spaces of type I curves at $\cS_{(\sdelta,\va,\vb)}$ span the
 component of $T_{(\sdelta,\va,\vb)}\cQ_d$ given by the first summand
 of~\eqref{Eq:Tangent_Space}. 
 (Recall that in type I, we have $\delta_i=\delta'_j=1$ and $\delta_j=\delta'_i=0$.)

 A similar analysis shows that the tangent space at $t=0$ of the type III curve
 defined by~\eqref{Eq:typeIII} is spanned by $x^{a'_i-a_i}y^{b'_i-b_i}\bE_{ii}$,
 and so the tangent spaces  at $\cS_{(\sdelta,\va,\vb)}$ of type III curves span the
 component of $T_{(\sdelta,\va,\vb)}\cQ_d$ given by the second summand
 of~\eqref{Eq:Tangent_Space} when $i=j$.

 For a curve of type II, note that the family of sheaves described
 by~\eqref{Eq:typeII} is constant in a neighborhood of $\infty$.
 In a neighborhood of $0$, it is given by
\[
   x^{a_i+c}\be_i,\quad x^{a_j}\be_j,\quad\mbox{ and }\quad 
   x^{a_i}\cdot\bigl(\be_i + t x^{(a_j-c)-a_i}\bE_{ij}(\be_i)\bigr)\,.
\]
 Thus $x^{(a_j-c)-a_i}\bE_{ij}$ spans the tangent space at $t=0$.
 A similar argument near $\infty$ for the sheaves described
 by~\eqref{Eq:typeII-b} shows that the tangent spaces of type II curves at
 $\cS_{(\sdelta,\va,\vb)}$ span the component of
 $T_{(\sdelta,\va,\vb)}\cQ_d$ given by the 
 second summand of~\eqref{Eq:Tangent_Space} when $i\neq j$.
\end{proof}


A {\em moment graph} of a $T$-variety is a graph 
whose vertices correspond to $T$-fixed points and whose edges correspond to
$T$-invariant curves, embedded into $\R\otimes\Hom(T,\Z)$ so that the edge
corresponding to a $T$-invariant curve is parallel to 
the weight of the action of $T$ on the curve.   
More specifically, if $C$ is a $T$-invariant curve joining fixed points $p$ and
$q$, then the edge from $p$ to $q$ in the moment graph is a positive 
multiple of the $T$-weight of $T_pC$.
When $k=\C$ and we fix a K\"ahler form, there is a moment map
$\mu\colon\cQ_d\to {\mathfrak t}^*$ and the image of the
$T$-fixed points and $T$-invariant curves is a moment graph.

When there are finitely many $T$-invariant curves, the Goresky-Kottwitz-MacPherson
method to compute equivariant cohomology is conveniently expressed in 
terms of a moment graph, with one relation for each edge.  
When there are infinitely many $T$-invariant curves, there are additional 
relations coming from families of $T$-invariant curves, so it is 
better to work with the moment {\em multigraph\/}, where each family of 
$T$-invariant curves (which will appear as a connected component of parallel
edges in the moment graph) is considered to form a single multiedge with
more than $2$ vertices, given by the fixed points in the closure of the 
family.  
To have a structure which determines the equivariant cohomology
or Chow groups, we should label each multiedge with the topological
type of the corresponding family.

Guillemin and Zara \cite{GZ01,GZ02,GZ03} have explored the combinatorial
properties of moment graphs.

\begin{example}\label{example_two}
{\rm 

Figure~\ref{F:Q202} represents a moment multigraph of $\cQ_2(0,2)$.
\begin{figure}[htb]
\[
  \begin{picture}(210,220)(0,-10)
   \put( 40, 10){\includegraphics[height=180pt]{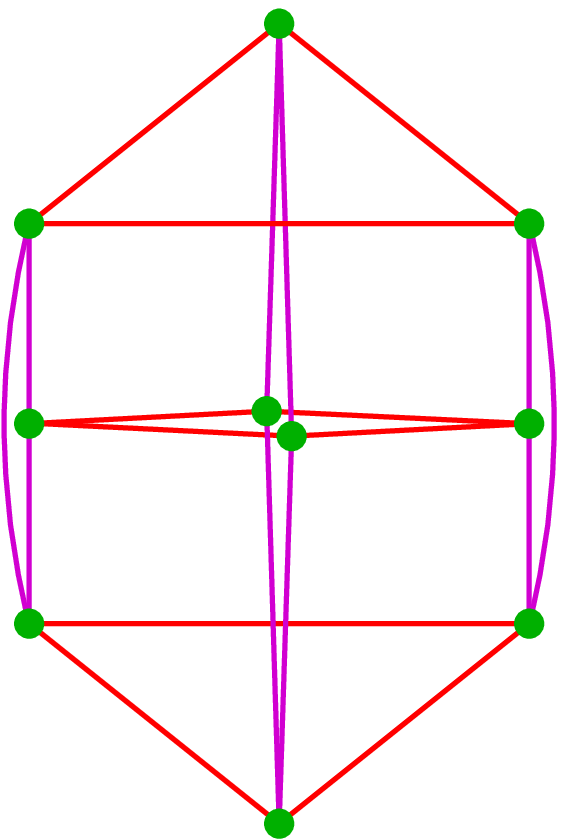}}
   \put( 10,140){\includegraphics[width=25pt]{figures/11.20.00.eps}}
   \put( 88,193){\includegraphics[width=25pt]{figures/11.11.00.eps}}
   \put(168,140){\includegraphics[width=25pt]{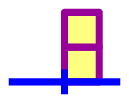}}
   \put( 10, 90){\includegraphics[width=25pt]{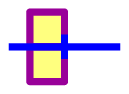}}
   \put(168, 90){\includegraphics[width=25pt]{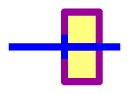}}
   \put( 10, 39){\includegraphics[width=25pt]{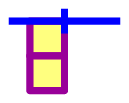}}
   \put( 88, -8){\includegraphics[width=25pt]{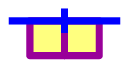}}
   \put(168, 39){\includegraphics[width=25pt]{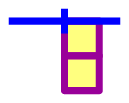}}
   \put(108, 78){\includegraphics[width=25pt]{figures/11.10.01.eps}}
   \put( 70,108){\includegraphics[width=25pt]{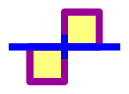}}

   \thicklines

   \put(-20,5){\vector(1,0){50}} \put(-12,-5){$\be_1{-}\be_2$}
   \put(-20,5){\vector(0,1){30}} \put(-30, 5){$\bff$}
  
  \end{picture}
\]
\caption{A moment multigraph of $\mathcal{Q}_2(0,2)$}
\label{F:Q202}
\end{figure}
Since the $T$-fixed points  
\raisebox{-5pt}{\includegraphics[width=20pt]{figures/11.10.01.eps}}
and 
\raisebox{-5pt}{\includegraphics[width=20pt]{figures/11.01.10.eps}}
have the same image in this and in any moment multigraph, 
we displace their images from their true positions for clarity. 
Similarly, some images of  $T$-invariant curves are displaced or
drawn as arcs.

The $T$-basis of the tangent space at  
\raisebox{-5pt}{\includegraphics[width=20pt]{figures/11.20.00.eps}}
has weights $\{ \be_1-\be_2, \be_1-\be_2+\bff, -\bff, -2\bff\}$.
These correspond to the following four $T$-invariant curves:
\[
  \begin{array}{|l|c|c|c|c|}\hline
   \mbox{Submodules of $\be_1\cO\oplus\be_2\cO$}&
    t=0&s=0&\mbox{Weight}&\mbox{Type}\\\hline 
    x^2\be_1, \ x^2\be_2,\ s\be_1+t\be_2& x^2\be_1,\ \be_2 
      \rule{0pt}{12pt}  & \be_1,\ x^2\be_2& \be_1-\be_2& {\rm II}\\\hline
    x^2\be_1,\ x^2\be_2,\ sx\be_1+ty\be_2&x^2\be_1,\ \be_2 & x\be_1,\ x\be_2
      \rule{0pt}{12pt}  & \be_1-\be_2+\bff&{\rm II}\\\hline
    (sx^2+txy)\be_1,\  \be_2 &x^2\be_1,\ \be_2 
      \rule{0pt}{12pt}  & xy\be_1 ,\ \be_2&-\bff&{\rm III}\\\hline
    (sx^2+ty^2)\be_1,\ \be_2 &x^2\be_1,\ \be_2 
     \rule{0pt}{12pt}   & y^2\be_1,\ \be_2&-2\bff&{\rm III}\\\hline
   \end{array}
\]
The $T$-basis to the tangent space at 
\raisebox{-5pt}{\includegraphics[width=20pt]{figures/11.10.01.eps}}
has weights $\{ \pm\bff, \pm(\be_2-\be_1)\}$.
These correspond to the following four $T$-invariant curves:
\[
  \begin{array}{|l|c|c|c|c|}\hline
   \mbox{Submodules of $\be_1\cO\oplus\be_2\cO$}&
    t=0&s=0&\mbox{Weight}&\mbox{Type}\\ \hline
    (sx+ty)\be_1,\ y\be_2 &x\be_1,\ y\be_2 & y\be_1,\ y\be_2
       &-\bff&{\rm III}\\\hline
    xy\be_1,\ y\be_2,\ sx\be_1+tx\be_2& x\be_1,\ y\be_2 
                       & xy\be_1,\ \be_2& \be_2-\be_1&{\rm II}\\\hline
    x\be_1,\ xy\be_2,\ ty\be_1+sy\be_2&x\be_1,\ y\be_2 & \be_1,\ xy\be_2
                       & \be_1-\be_2&{\rm II}\\\hline
    x\be_1,\ (sy+tx)\be_2 &x\be_1,\ y\be_2 & x\be_1,\ x\be_2
       &\bff&{\rm III}\\\hline
   \end{array}
\] 
 
}\end{example}

\section{Families of $T$-invariant curves}
\label{Tfamilies}

Suppose that $Z$ is a $T$-invariant curve on $\cQ_d$.
Let $T'$ be the identity component of the stabilizer of $Z$ in $T$, a
codimension one subtorus of $T$. 
Then the action of $T$ on $Z$ factors through the quotient
\[
   T\ \longrightarrow\ T/T'\ \simeq\ \G_m\,.
\]
This composition $\chi\colon T\to \G_m$ is a primitive weight parallel to the
weight of the action of $T$ on $Z$.
Let $Y$ be the component of the $T'$-fixed point locus $\cQ_d^{T'}$ which
contains $Z$.
Then $Y\setminus Y^T$ is foliated by one-dimensional orbits of $T$
whose closures are $T$-invariant curves.
We call $Y$ the family of $T$-invariant curves on $\cQ_d$ which contains $Z$.
If $p\in Y^T\subset\cQ_d^T$ is a $T$-fixed point of $Y$, then $T_pY$ is a
$T$-invariant linear subspace of $T_p\cQ_d$ which is fixed pointwise by $T'$.
In particular, all weights of the $T$-action on $T_pY$ are multiples of $\chi$. 

To classify families of $T$-invariant curves, we first 
determine which $T$-weights of $T_p\cQ_d$ are parallel.
If the weight of a curve is not parallel to the weight of any other curve, it is
isolated. 
By Theorem~\ref{tangentdirections}, it suffices to identify those $T$-invariant
curves incident on a fixed point with parallel $T$-weights. 

\begin{thm}\label{T:parallel}
  Two $T$-invariant curves $\cS$, $\cS'$ of types {\rm I}, {\rm II}, 
  or\/ {\rm III} containing the fixed point $\cS_{(\sdelta,\va,\vb)}$ have
  parallel $T$-weights if and only if either 
 \begin{enumerate}
  \item both $\cS$ and $\cS'$ have type {\rm III} or
  \item both $\cS$ and $\cS'$ have type {\rm II}, and 
   \begin{enumerate}
    \item $\cS$  connects $\cS_{(\sdelta,\va,\vb)}$ to $\cS_{(\sdelta,\va',\vb)}$,
    \item $\cS'$ connects $\cS_{(\sdelta,\va,\vb)}$ to $\cS_{(\sdelta,\va,\vb')}$,
    \item $\va$ and $\va'$ agree except in positions $i$ and $j$ with 
$i\neq j$; $\vb$ and $\vb'$  agree except in positions $i$ and $j$ (same $i,j$);
 and $a_i+b_i=a'_j+b'_j$.
   \end{enumerate}
 \end{enumerate}
\end{thm}

\begin{proof}
 For (1), note that the weight  of $T$ on a type III curve is parallel to
 $\bff$.
 
 If a curve does not have type III, then its weight has the form
 $\be_j-\be_i+c\,\bff$, where $\delta_i=1$ and either $\delta_j=0$ if it has type
 I or $\delta_j=1$ if it has type II.
 Thus $\cS$ and $\cS'$ have the same weight and type, $\vdelta=\vdelta'$, and
 the indices $i$ and $j$ in their definitions coincide.

 Weights of curves of types I and II correspond to~\eqref{Eq:TI}
 and~\eqref{Eq:TII_III}, respectively. 
 Inspecting~\eqref{Eq:TI}, shows that no two curves of type I can have the same weight.
 Inspecting~\eqref{Eq:TII_III} reveals that either a given curve $\cS$ of type
 II has a unique weight, or else there is exactly one other type II curve $\cS'$ with
 the same weight, and the two curves are as described in the statement of the
 theorem. 
\end{proof}

We show that all type III curves at a fixed point lie in a single family of
$T$-invariant curves, and if two type II curves
have the same weight then they lie in a 2-dimensional family.
Together with the isolated $T$-invariant curves,
this shows that the tangent spaces at a given fixed point
$T$ to families of $T$-invariant curves are the subspaces of $T_p\cQ_d$
which are stabilized by codimension 1 subtori of $T$.
It follows that these families contain all $T$-invariant curves in $\cQ_d$.
\medskip

\noindent{\bf Vertical families.}  
A \emph{vertical family} containing $\cS_{(\sdelta,\va,\vb)}$ is parametrized by
the product of projective spaces 
\[
    \prod_{\delta_i=1} \P H^0(\cO(a_i+b_i))\ \simeq\ 
    \prod_{\delta_i=1} \P^{a_i+b_i}\,.
\]  
It contains exactly the fixed points $\cS_{(\sdelta',\va',\vb')}$ where
$\vdelta=\vdelta'$, and $\va'+\vb'=\va+\vb$, along with all
type III curves which connect them.  
These include all type III curves at each of these fixed points.

Consider the family $\cS$ of submodules of $\cO^n_{\P^1}$ generated by
\[
   \{  \be_i s_i  \mid \delta_i=1,\quad s_i\in H^0(\cO(a_i+b_i))\}\,.
\]
The base of this family is $\prod_{\delta_i=1} \P H^0(\cO(a_i+b_i))$,
all subsheaves have rank $n-r$ and degree $-d$, and the foliation by
$T$-invariant curves is given by the $T$-action on the base.
\medskip

\noindent{\bf Horizontal families.}
If there exist $i,j,c,c'$ such that $1 \le c \le a_j$, 
$1 \le c' \le b_j$, and $a_i+b_i+c+c'=a_j+b_j$, then the point
$\cS_{(\sdelta,\va,\vb)}$ 
lies in a \emph{horizontal family} parametrized by the product of two
projective lines.  Let $a'_i=a_i+c$, $a'_j=a_j-c$, $b'_i=b_i+c'$, 
$b'_j=b_j-c'$.
Let $([s,t], [\sigma,\tau])$ be the coordinates of
$\P^1\times\P^1$, and let $\cS$ be the submodule of $\cO^n_{\P^1}$
which, except for its components in $\cO\cdot\be_i + \cO\cdot\be_j$,
agrees with $\cS_{(\sdelta,\va,\vb)}$.
The component of $\cS$ in $\cO\cdot\be_i + \cO\cdot\be_j$ is the subsheaf
generated by 
 \[
   \be_jx^{a_j}y^{b_j}, \ \  
   s\be_ix^{a'_i}y^{b_i}+t\be_jx^{a_j}y^{b'_j},\ \ 
   \sigma \be_ix^{a_i}y^{b'_i}+\tau\be_jx^{a'_j}y^{b_j},\ \ 
   \be_ix^{a'_i}y^{b'_i}\,.
 \]
Similar reasoning as for Theorem~\ref{Th:Curves_are_correct} shows that 
this defines a family of $T$-invariant
curves over the base $\P^1\times\P^1$ with coordinates $([s,t], [\sigma,\tau])$.  
It contains
four $T$-fixed points:  setting $t=\tau=0$ gives the fixed point
$\cS_{(\sdelta,\va,\vb)}$, setting $t=\sigma=0$ gives the fixed point
$\cS_{(\sdelta,\va',\vb)}$, setting  $s=\tau=0$ gives the fixed point
$\cS_{(\sdelta,\va,\vb')}$, and  setting $s=\sigma=0$ gives the fixed point
$\cS_{(\sdelta,\va',\vb')}$.  This family also contains 
the four Type II curves connecting these four fixed points, given by setting 
exactly one of $s,t,\sigma$, or $\tau$ equal to zero.
Furthermore, the data $(\vdelta,\va,\vb)$ and $(\vdelta',\va',\vb')$ satisfy 
Theorem~\ref{T:parallel}(2)(c), and any two $T$-invariant curves $\cS$ and
$\cS'$ as in Theorem~\ref{T:parallel}(2) lie in a unique horizontal family.


\section{Algebraic extension of GKM theory}\label{extended GKM theory}
We discuss equivariant localization and an extension
of the Goresky-Kottwitz-Mac\-Pherson relations 
when there are finitely many fixed points
but infinitely many $T$-invariant curves.  We work with 
equivariant Chow rings; similar results hold for equivariant
cohomology.
In fact, when $k=\C$, $X$ is smooth and projective, and $X^T$ is finite,
the two theories coincide.

We first recall some properties of $T$-equivariant Chow rings as
developed by Edidin and Graham~\cite{EG98a} and Brion~\cite{brion}.
Next, we outline Evain's~\cite{evain} development of ideas of Brion which 
extends the GKM relations to 
describe the $T$-equivariant Chow ring of a smooth variety with finitely
many $T$-fixed points when there are infinitely many $T$-invariant curves.
This description involves ideal-membership relations, one for each generator of
the equivariant Chow ring of each family of $T$-invariant curves.
When the generators are given by smooth subvarieties, these 
relations may be expressed in terms of tangent weights.
This gives one form of our presentation for $A^*_T(\cQ_d)$ in
Theorem~\ref{Main_Theorem}. 
We next give a variant of these relations using 
differential operators, which gives the other form  of our
presentation for $A^*_T(\cQ_d)$.  
We then compute these relations for products of projective spaces, 
and finally deduce Theorem~\ref{Main_Theorem}.

\subsection{Torus equivariant Chow rings}

When a linear algebraic group $G$ acts on a smooth scheme $X$, 
Edidin and Graham~\cite{EG98a} defined the equivariant Chow ring
$A^*_G(X)$, using Totaro's algebraic approximation to the classifying space 
of $G$.   It satisfies 
functorial properties under equivariant maps analogous to 
those for ordinary Chow rings~\cite{Fu98}, including proper pushforwards
and pullbacks by local complete intersection morphisms.

When the group is a torus $T$, Brion \cite{brion} gave an alternative 
development of this theory which includes versions of the localization
theorems that hold for equivariant cohomology. 
He gave the following presentation for the equivariant Chow ring, 
analogous to the usual presentation of Chow groups.
The equivariant Chow ring $A^*_T(p)$ of a point $p$ is the integral symmetric
algebra $S$ of the character group $\hat{T}$ of $T$.
Equivariant pullback makes $A^*_T(X)$ into an $S$-module.

\begin{prop}[Theorem~2.1~\cite{brion}]\label{P:Chow_Gens}
 The $S$-module $A^*_T(X)$ is defined by generators $[Y]$, for each 
 $T$-invariant subvariety $Y$, and by relations 
 $[\mbox{div}_Y(f)]-\chi[Y]$  for each rational function $f$ on $Y$ which is
 a $T$-eigenvector of weight $\chi$; here $\chi$ is considered as
 an element of $S$ in degree $1$.
\end{prop}

It follows immediately that the usual Chow ring may be recovered 
from the $S$-module $A^*_T(X)$, as the quotient
by the ideal $S^+$ of $S$ generated by the character group $\hat{T}$.

\begin{prop}[Corollary~2.3.1~\cite{brion}]\label{P:Chow}
  We have
\[
   A^*(X)\ =\ A^*_T(X)\otimes_S \Z\ =\ A^*_T(X)/S^+ A^*_T(X)\,.
\]
\end{prop}

The analogous statement in equivariant cohomology requires stronger 
hypotheses. 

When $k=\C$,  the connection between 
Chow groups and cohomology is given by the cycle map
 \[
   A^*_T(X)\ \to\  H^*_{T,c}(X, \Z)
 \]
to compactly supported (Borel-Moore) equivariant cohomology.
If $X$ is projective and the fixed point set $X^T$ is finite, then
the cycle map is an isomorphism.

Some statements below hold only for the rational
equivariant Chow ring $A^*_T(X)_\Q := A^*_T(X) \otimes_\Z \Q$.
This is a module over the rational equivariant Chow ring 
$A^*_T(p)_\Q$ of a point $p$, which is the symmetric algebra $S_\Q$ of
$\hat{T}_\Q := \hat{T}\otimes_\Z \Q$.  

\subsection{Localization}
We now assume that $X^T$ is finite, and that $X$ has a 
decomposition into $T$-invariant affine cells $C_1,\dotsc,C_m$ which can be
ordered so that for $i=1,\dotsc,m$, the union $C_1\cup\dotsb\cup C_i$ 
is Zariski open.
We will call such varieties {\em filtrable}; this is close to the terminology 
Brion used in \cite{brion}, but he did not require that
$X^T$ be finite, and his cells were allowed to be vector bundles over components
of $X^T$. 
If $X$ is smooth and projective and $X^T$ is finite, then 
Bia{\l}ynicki-Birula~\cite{BB73} showed that it is filtrable.

Let $i\colon X^T\to X$ be the inclusion of the subscheme of
$T$-fixed points of $X$.  
\begin{prop}[Corollary 3.2.1~\cite{brion}] 
The $S$-module $A^*_T(X)$ is free.  The map
\[i^*\colon A^*_T(X) \to A^*_T(X^T)\]
is an injection.  
\end{prop}

Brion also established Chow ring versions of results of Chang and Skjelbred
and of Goresky, Kottwitz, and MacPherson concerning the image of the
localization map.
\begin{prop}[Sections 3.3 and 3.4~\cite{brion}]\label{P:GKM_algebraic}
\mbox{\ }
 \begin{enumerate}
  \item[(a)]
   The image of the localization map $i^* A^*_T(X)_\Q \to A^*_T(X^T)_\Q$ is the
   intersection of the images of the localization maps
\[
   i^*_{T'}\ \colon\ A^*_T(X^{T'})_\Q\ \to\  A^*_T(X^T)_\Q
\]
  where $T'$ runs over all codimension one subtori of $T$.

 \item[(b)]
  When $T$ acts with finitely many fixed points and has finitely
  many invariant curves, then the image of the localization map 
  $i^*\colon A^*_T(X)_\Q \to A^*_T(X^T)_\Q\simeq (S_\Q)^{X^T}$
  is the set of all tuples $(f_p)_{p\in X^T} \in (S_\Q)^{X^T}$ such that
  whenever $p$ and $q$ belong to the
  same irreducible $T$-invariant curve $C$, we have 
  $f_p\equiv f_q$ modulo $\chi$, where $\chi$ is the weight of
  the action of $T$ on $T_pC$.
 \end{enumerate}
\end{prop}

Statement (a) is analogous to a theorem of 
Chang and Skjelbred \cite{ChSk74} for equivariant cohomology.
This result, together with the easy calculation of the equivariant
Chow groups of $\P^1$, immediately gives (b), which is the Chow analog of 
the GKM relations for equivariant cohomology.

In general this result does not hold with $\Z$ coefficients.  For instance, 
suppose that $\dim X = 2$,  $x\in X^T$, and the weights of $T$ on the
tangent space $T_xX$ are $a\chi$ and $a'\chi'$, where 
$\chi, \chi'\in \hat{T}$ are linearly independent primitive characters.  
Then condition (b) would say that if $(f_p) \in A^*_T(X^T)$ 
has $f_p = 0$ for $p \ne x$, then it is in the image of $i^*$ if $f_x$ is a 
multiple of $\lcm(a,a')\chi\chi'$.  In fact, $f_x$ must
be a multiple of $aa'\chi\chi'$.

This is essentially the only obstruction to working with $\Z$ coefficients,
at least if the fixed point set is finite.  We say that the tangent weights
of a $T$-variety $X$ are {\em almost coprime} if whenever two 
$T$-weights of $T_pX$ for $p\in X^T$ are divisible by the same 
integer $a > 1$, then they are parallel. 
%
%
With this added hypothesis, Brion's proof
of Proposition~\ref{P:GKM_algebraic} given by Brion works over $\Z$.

\begin{thm} \label{T:GKM_algebraic_Z}
 Let $X$ be a smooth filtrable $T$-variety whose 
 tangent weights are almost coprime.
 Then Proposition~$\ref{P:GKM_algebraic}$ holds with rational Chow groups
 replaced by integral Chow groups.
\end{thm}

\subsection{Evain's relations}
When $T$ does not have finitely many invariant curves on $X$, then
statement (b) of Proposition~\ref{P:GKM_algebraic} fails, but by (a) we
can still compute $A^*_T(X)$ if we know the images of $i^*_{T'}$ for 
all codimension one subtori $T'$ of $T$.  A finite set of such $T'$
suffices, namely those which fix at least one $T$-invariant curve
pointwise.  Brion \cite{brion} and Goldin and Holm \cite{GH} have 
computed cases where the components of $X^{T'}$ are low-dimensional.
Evain \cite{evain} recently described relations in the general case:
we recall his results.

Let $Y = X^{T'}$ for $T'$ a codimension one subtorus of $T$.
By \cite{I}, $Y$ is smooth.
For $p\in Y^T = X^T$, let $e^T_p(Y) = e^T(TY)|_p$ be the localization of the
equivariant Euler class of $TY$ at $p$.  
Under the identification $A^*_T(p) = S$, this 
is the product of the $T$-weights on the tangent space $T_pY$.

\begin{prop}[Corollary 27~\cite{evain}]\label{P:C27Ev}
 A class $\alpha=(\alpha_p)_{p\in Y^T}$ in $S^{Y^T}$ 
 lies in $i^*_{T'}A^*_T(Y)$ if and only if 
\begin{equation}\label{E:Evain}
  \sum_{p\in Y^T} \frac{\alpha_p\beta_p}{e^T_p(Y)}\ \ \in\ S
\end{equation}
 for every $\beta\in i^*_{T'}A^*_T(Y)$.
\end{prop}

\noindent{\bf Remark on Evain's proof.}
The condition~\eqref{E:Evain} is necessary, since
if $\pi$ is the 
projection of $Y$ to a point,
 then  
 the sum is simply $\pi_*(\alpha\cdot\beta)$, by the integration
 formula of Edidin and Graham~\cite{EG98b}.
Note that since $\pi_*$ is $S$-linear, it is enough to take
$\beta$ in a generating set of the $S$-module
$i^*A^*_T(Y)$.

 By Bia{\l}ynicki-Birula \cite{BB73}, there are two
 $T$-invariant cell decompositions $C_p^+$ and $C_p^-$ for $p\in Y^T$ of $X$ and
 an ordering of the fixed points $Y^T$ such that the matrix with entries in $S$
 whose $(p,q)$-entry is 
\[
   \pi_*([C^+_p]\cdot [C^-_q])
\]
 is unitriangular.
 Either set of classes $[C^+_p]$ or $[C^-_p]$ forms a basis for the
 $S$-module $A^*_T(X)$, and
 expressing the elements $\alpha$ and $\beta$ in these two bases proves
 sufficiency.  \medskip

Combining Proposition~\ref{P:C27Ev} with Proposition \ref{P:GKM_algebraic}
and Theorem~\ref{T:GKM_algebraic_Z} gives the
following criterion for membership in $i^*A^*_T(X)$.

\begin{thm}\label{T:E_Criterion}
  A class $\alpha=(\alpha_p)_{p\in X^T} \in (S_\Q)^{X^T}$ lies in the image
  $i^*A^*_T(X)_\Q$ of the localization map if and only if
for all $Y = X^{T'}$ for $T'$ a codimension one subtorus of $T$ we have
\[
   \sum_{p\in Y^T} \frac{\alpha_p\beta_p}{e^T_p(Y)}\ \ \in\ S_\Q\,
\]
for all $\beta$ in a set of $S_\Q$-module generators 
for $i^*_{T'}A^*_T(Y)_\Q$.


If the tangent weights of $X$
are almost coprime, the same statement holds over $\Z$.
\end{thm}

\begin{rmk}
  When $X$ is smooth, the relations in
  Theorem~\ref{T:E_Criterion} can also be taken for $Y$ running over
  all irreducible components of the union of the fixed points
  and the $T$-invariant curves, since such $Y$ are just the connected
  components of the $T'$-fixed loci $X^{T'}$ for some codimension one subtorus
  $T'$ of $T$.  
  We call this union of fixed points and $T$-invariant curves the 
  {\it one-skeleton of $X$}. 
\end{rmk}

To apply Theorem~\ref{T:E_Criterion}, we need to know
explicit generators of $A^*_T(Y)$, or more precisely their localizations
to $Y^T$.
By Proposition~\ref{P:Chow_Gens}, one class of generators are the equivariant
fundamental cycles $[Z]$ of $T$-invariant subvarieties $Z$ of the components
$Y$.  These are easy to compute when $Z$ is smooth, since if $p\in Z^T$ we have 
$[Z]_p = e^T_p(N_ZY)$, the equivariant
Euler class of the normal bundle to $Z$ in $Y$, while
if $p\in Y^T \setminus Z^T$, then $[Z]_p = 0$.
It follows that $\frac{[Z]_p}{e^T_p(Y)}=\frac{1}{e^T_p(Z)}$ if $p\in Z^T$.

%
%
To see this, note that $[Z]_p\in A^*_T(p)$ is the pullback of $[Z]$ along
the regular embedding $i_{p,Y}\colon p\to Y$.
We factor $i_{p,Y}$ as the composition 
\[
   p\ \xrightarrow{\ i_{p,Z}\ }\ Z
      \xrightarrow{\ i_Z\ }\ Y\,.
\]
The class $[Z]\in A^*_T(Y)$ is the pushforward along $i_Z$ of the unit
class $1=[Z]\in A^*_T(Z)$, and so we have
\[
  [Z]_p\ =\ i^*_{p,Z} i^*_Z i_{Z,*} 1\ =\ 
  i^*_{p,Z} e^T(N_ZY) =\ e^T_p(N_ZY)\,,
\]
by the self-intersection formula for Chow rings. 

Thus if we can find for each $Y$ a collection $\cZ_Y$ of smooth $T$-invariant
subvarieties of $Y$ so that the classes $[Z]$ for $Z\in \cZ_Y$ generate
$A^*_T(Y)$ as an $S$-module, we get the following more explicit version of
Theorem \ref{T:E_Criterion}. 

\begin{thm}\label{T:toric_criterion}
 A class $\alpha=(\alpha_p)_{p\in X^T}\in (S_\Q)^{X^T}$ lies in
 $i^*A^*_T(X)_\Q$ if and only if 
 \begin{equation} \label{relation for smooth Z}
   \sum_{p\in Z^T} \frac{\alpha_p}{e^T_p(Z)}\ \in\ S_\Q\,
 \end{equation}
 for all $Z \in \cZ_Y$ and all components $Y$ of the one-skeleton of $X$. 
 If the tangent weights of $X$ are almost coprime, the same statement holds over $\Z$. 
\end{thm}

 The necessity of \eqref{relation for smooth Z} 
 does not require the argument above,
 since if $(\alpha_p) = i_*\alpha$ for $\alpha \in A^*_T(X)$, then the sum
 is just $\pi_*(\alpha|_Z)$, where $\pi$ is the projection of $Z$ to a
 point.

Obvious candidates for the subvarieties $[Z]$ are the
closures of the Bia{\l}ynicki-Birula cells, since their classes
form  an $S$-basis for $A^*_T(Y)$.  Unfortunately, they are not  
in general smooth---this was the case for Evain.
However, for the quot schemes we study, they are smooth, as the 
connected components $Y$ are products of projective spaces.  
More generally we can ask that for each component $Y$ of 
$X^{T'}$ there is a torus $T_Y$ containing $T$ which acts on $Y$ 
with finitely many orbits, so that $Y$ is a smooth toric variety. 
The closures of the cells will be $T_Y$-orbit closures, and therefore
smooth.

The relations of Theorem \ref{T:toric_criterion} 
are the same as those found by Goldin and Holm \cite{GH}
for equivariant cohomology of Hamiltonian $T$-spaces.
Their result applied where the spaces $X^{T'}$ 
are at most four-dimensional (over $\R$).  Since this
implies that the components of $X^{T'}$ are toric manifolds, 
Theorem~\ref{T:toric_criterion} 
recovers their result. 

\subsection{Evain's relations as differential operators}
We rewrite this algebraic criterion in a different form.
Suppose that $Y$ is a smooth component of $X^{T'}$ and $Z\subset Y$ is a smooth
$T$-invariant subvariety.  
The action of $T$ on $Y$ factors through a character $\eta\colon T \to \C^*$, 
so the weights of
$T$ on $T_pZ$ for $p\in Z^T$ are non-zero scalar multiples of $\eta$.
Thus there exist numbers $d_p=d_p(Z)$ so that 
\[
   e^T_p(Z) = d_p(Z)\cdot \eta^{\dim Z}\,.
\]
The terms in \eqref{relation for smooth Z} have a common
denominator $\eta^{\dim Z}$, and so we may rewrite it
as
\[
    \sum_{p\in Z^T} \frac{\alpha_p}{d_p(Z)}\ \in\ {\eta^{\dim Z}}S_\Q\,.
\] 

We can rewrite this condition using a linear differential operator. 
The ring $S_\Q$ is the symmetric algebra of $\hat{T}_\Q$, or dually
the ring of polynomial functions on ${\hat{T}_\Q}^*$.
Choose $\zeta\in {\hat{T}_\Q}^*$
for which $\zeta(\eta)\neq 0$.  Then 
the operator $D=D_\zeta$ of differentiation in the direction of
$\zeta$ acts on $S_\Q$.  If $f\in S_\Q$ is divisible by $\eta$, then
$\eta^k$ divides $f$ if and only if $\eta^{k-1}$ divides $Df$, so
the relation \eqref{relation for smooth Z} is equivalent to
 \[ 
   \sum_{p\in Z^T} d_{p}(Z)^{-1}D^j \alpha_p \equiv 0
   \mod\eta,
    \quad\mbox{for all}\ 0 \le j < \dim Z.
\]

We give a variant of Theorem~\ref{T:toric_criterion} which uses the
last relation, but only with the maximum order derivative $j = \dim Z - 1$.
In exchange, we must apply it using more subvarieties $Z$. 

Let $Y$ be a component of the one-skeleton of $X$, let 
$\eta$ be the associated character of $T$, and consider
the two Bia{\l}ynicki-Birula cell decompositions
$\{C^-_p\mid p\in Y^T\}$ and $\{C^+_p\mid p\in Y^T\}$ 
induced by the $T$-action.
Each cell $C=C^+_p, C^-_p$ is isomorphic to the $T$-vector space 
$T_pC\subset T_pY$.
Suppose that within each cell $C=C^+_p$ 
we can find $T$-invariant affine subspaces $C_{p,1}, \dots C_{p,\dim C}$ with
$\dim C_{p,i} = i$ and which have smooth closures $Z_{p,i} = \overline{C_{p,i}}$.
As before, this will be true if each $Y$ is a toric variety for a larger torus 
$T_Y$ containing $T$, since we can take each $Z_{p,i}$ to be the closure of a $T_Y$-orbit.

\begin{thm} \label{derivative theorem}
With these assumptions, a 
class $\alpha=(\alpha_p)_{p\in X^T}\in (S_\Q)^{X^T}$ lies in
$i^*A^*_T(X)_\Q$ if and only if 
 \begin{equation} \label{top degree derivative relation}
  \sum_{q\in Z^T} d_{q}(Z)^{-1}D^{\dim Z - 1} \alpha_q \equiv 0 \mod\eta .
 \end{equation}
for all $Z = Z_{p,i}$ and for all components $Y$ of the one-skeleton of $X$.
\end{thm}

\begin{proof}
The necessity of the conditions \eqref{top degree derivative relation} follows from the 
previous discussion.  

To show they are sufficient, let $U$ be an open union of the 
%
cells $C^-_p$ and  note that
$Z_{p,i}\subset U$ if and only if $p \in U$.  
We use induction on the number of cells in $U$ to show that the image of 
$i_U^*\colon A^*_T(U)_\Q \to A^*_T(U^T)_\Q$ is the set of $(\alpha_x)|_{x\in U^T}$
satisfying \eqref{top degree derivative relation} for all $Z_{p,i}\subset U$.

When $U$ is a single cell, this is immediate, as $i^*_U$ is an isomorphism
and there are no $Z_{p,i}$'s contained in $U$.  
Otherwise, suppose
$\alpha = (\alpha_x)|_{x\in U^T}$ satisfies \eqref{top degree derivative relation}
for all $Z_{p,i}\subset U$.
Let $C^-_p \subset U$ be a closed cell, and put $U' = U \setminus C^-_p$.
%
%
There is an exact sequence\cite[Proposition~3.2]{brion} 
 \begin{equation} \label{exact sequence} 
  0\ \to\ A^*_T(U)_\Q\  
  \xrightarrow{\, \rho \,}\ A^*_T(U')_\Q \times A^*_T(C^-_p)_\Q \ 
   \to\ A^*_T(C^-_p)_\Q/(e^T(N))\ \to\ 0\,,
 \end{equation}
where $e^T(N)$ is the equivariant Euler class of the normal bundle $N$
to $C^-_p$ in $X$.  Under the isomorphism $A^*_T(C^-_p) \cong S$, this is just
the product of all the $T$-weights of $N$.  
The components of $\rho$ are the restriction maps, while the map
$A^*_T(C^-_p)_\Q \to A^*_T(C^-_p)_\Q/(e^T(N))$ is the natural quotient.

By the inductive hypothesis, $\alpha|_{(U')^T}$ lies in the image of $i^*_{U'}$.
The exact
sequence~\eqref{exact sequence}  implies that $A^*_T(U) \to A^*_T(U')$
is surjective  and so we can write $\alpha = i^*_U\beta + \gamma$, with
$\beta \in A^*_T(U)$ and $\gamma|_{(U')^T} = 0$. 
Since $i^*_U\beta$ satisfies the relations 
\eqref{top degree derivative relation}, so does $\gamma$.
It will be enough to show that $\gamma$ is in the image of $i^*_U$.  But using 
the exact sequence \eqref{exact sequence}, we see that this holds 
if and only if $\gamma_p$ is a multiple of $c_d^T(N)$, which is a non-zero multiple
of $\eta^d$, where $d = \mathop{\mathrm{codim}} C^-_p = \dim C^+_p$.  But 
the relation \eqref{top degree derivative relation} implies that 
$D^{k}\gamma_p \equiv 0 \pmod \eta$ for $0 \le k < d$.  The result follows. 
\end{proof}


\begin{example}
\label{Example_1}
{\rm 
 Let $T = \G_m$ act on $X = \P^r$ by 
 \[
   t\cdot[x_0:x_1:\dots:x_r]\ =\ [x_0: tx_1: \dots: t^rx_r]
 \]
in homogeneous coordinates, where $t\in T$.
For each $0\le j \le r$, let $p_j \in X^T$ denote the $T$-fixed
point corresponding to the $j$th standard basis vector $\be_j$.  The
tangent space $T_{p_j}X$ is $\C^n/\C\cdot\be_j$, with the action of
$T$ given by $t\cdot \overline{\be}_k = t^{k-j}\overline{\be}_k$.
Thus $e_{p_j}(X) = (-1)^j j!(r{-}j)!\eta^r$, where 
$\eta$ is the identity character.  
More generally, if $0\le l\le n \le r$, let 
$Z_{l,m} = \P\Span\{\be_l, \be_{l+1},\dots,\be_m\}$.  The same 
calculation gives 
 \[
    e_{p_j}(Z_{l,m})\ =\ (-1)^{j-l} (j-l)!(m-j)!\eta^{m-l}\,.
\]

We can apply Theorem \ref{T:toric_criterion} using the smooth
subvarieties $Z_{0,l}$, for $1\le l \le r$.
Then $i^*A^*_T(X) \subset A^*_T(X^T)$
is the set of tuples $\alpha = (\alpha_0,\dots,\alpha_r)$ where
 \begin{equation}\label{example relation 1} 
     \sum_{0\le j\le l} \frac{(-1)^j\alpha_j}{j!(l-j)!}\ \in\ \eta^l S\,
 \end{equation}
for all $1 \le l \le r$.

On the other hand, we can apply Theorem \ref{derivative theorem}
using all the subvarieties $Z_{l,m}$. 
If $D$ is differentiation on $S$ in the direction of $\eta^\vee$, then
$\alpha$ lies in the image of $i^*$ if and only if 
 \begin{equation}\label{example relation 2}
   \sum_{l\le j\le m}
   \frac{(-1)^{j-l}D^{m-l-1}\alpha_j}{(j-l)!(m-j)!}
     \ \in\ \eta S
 \end{equation}
for all $0 \le l < m \le r$.  
We could take one more derivative and
ask that the resulting sums vanish, but this would not generalize
to actions of higher dimensional tori.

When $X = \P^1$ we get exactly the GKM relation for a primitive action.
}
\end{example}

\begin{example}\label{Example_2}
{\rm 
 These same arguments apply to products of projective spaces.   
Let $\vr = (r_1,\dots,r_n)$, and let $X = \P^{r_1} \times \dots \times 
\P^{r_n}$, where the action of $t \in T = \G_m$ on a point
$([x^1_0:x^1_1:\dots:x^1_{r_1}], \dots,[x^n_0:x^n_1:\dots:x^n_{r_n}])$ is
given by multiplying $x^i_j$ by $t^j$.  The fixed points have
the form $p_{\vj} = (p_{j_1},\dots,p_{j_n})$, where $p_{j_i}$ is the
$j_i$th fixed point in $\P^{r_i}$, in the notation of Example~\ref{Example_1},
and $\vj = (j_1,\dots,j_n)$ satisfies $0\le \vj \le \vr$, meaning that
$0\le j_i \le r_i$ for all $1 \le i \le n$.

For $\vl,\vm \in \Z^n$ with $0 \le l_i \le m_i \le r_i$ for all $i$, set
$Z_{\vl,\vm} = Z_{l_1,m_1}\times\dots \times Z_{l_n,m_n}$.  
For each $\vl \le \vj \le \vm$, the tangent space to $Z_{\vl,\vm}$ at
$p_{\vj}$ is $\bigoplus_{i = 1}^n T_{p_{j_i}}Z_{l_i,m_i}$.  
Using the computation from Example~\ref{Example_1},
we see that 
 \[
   e^T_{p_{\vj}}(Z_{\vl,\vm})\ =\
   (-1)^{|\vj|}(\vj-\vl)!(\vm-\vj)!\eta^{|\vm|-|\vl|}\,.
 \]
Recall that 
for an $n$-tuple $\va = (a_1,\dots,a_n)$, we put $\va! = a_1!\cdots a_n!$.

As in Example~\ref{Example_1}, we can either apply 
Theorem \ref{T:toric_criterion} using the 
subvarieties $Z_{\vzero,\vl}$, or Theorem \ref{derivative theorem} using 
all the $Z_{\vl,\vm}$.  The resulting conditions for a tuple 
$\alpha = (\alpha_\vj)_{\vzero\le \vj\le \vr}$
to be in the image of the localization map are just 
\eqref{example relation 1} and \eqref{example relation 2}, 
where the variables now represent elements of $\Z^n$ rather than scalars.
}
\end{example}

\subsection{Proof of Theorem \ref{Main_Theorem}}\label{main proof}
We combine these localization results with the 
geometry of the quot scheme from Sections
\ref{fixed_points}, \ref{Tcurves}, and \ref{Tfamilies} 
to produce a proof of Theorem \ref{Main_Theorem}.

The only non-primitive tangent weights are those with 
$i=j$ in~\eqref{Eq:TII_III}, which are multiples of $\bff$;
These correspond to Type III curves.  Thus the 
tangent weights are almost coprime, so 
Theorem \ref{T:toric_criterion} gives a
correct description of the integral equivariant Chow ring.

The relations I, II(a) and II(b) are the GKM relations
for the $T$-invariant curves of types I and II, as described
in Section \ref{Tcurves} using the identification
of the $T$-weights of these curves with the tangent weights
given by Theorem \ref{tangentdirections}.

The relations II(c)/II(c)$'$ come from the horizontal families of 
$T$-invariant curves of type II.  As described in 
Section~\ref{Tfamilies}, these are isomorphic to $\P^1 \times \P^1$,
where the action on 
each factor is by the same 
primitive character.  As in Example~\ref{Example_2},
we can apply the relation in Theorem~\ref{T:toric_criterion}
with $Y = \P^1 \times \P^1$ to get the relations II(c)/II(c)$'$.  
We get one new relation using $Z = Y$; smaller
$T$-invariant subvarieties contained in $Y$ are either $T$-invariant
curves, whose relations are already covered by II(a) and II(b), 
or points, which give no relation.  

Finally, the relations III and III$'$ come from the vertical families.  
As described in Section \ref{Tfamilies}, the vertical 
family containing $\cS_{(\sdelta,\va,\vb)}$ is isomorphic
to 
\[
    \prod_{\delta_i=1} \P H^0(\cO(a_i+b_i))\ \simeq\ 
    \prod_{\delta_i=1} \P^{a_i+b_i}\,,
\]   
and the fixed points in the family are those $\cS_{(\sdelta,\va',\vb')}$
with $\va' + \vb' = \va + \vb$.  
The codimension one subtorus $T_{k^n}\subset T$
acts trivially on the family, and the remaining action of 
$T_{\P^1}$ is the one described in Example~\ref{Example_2}, using the
monomial basis of $H^0(\cO(a_i+b_i))$. 
Example~\ref{Example_2} then gives exactly the relations III and III$'$.
This proves Theorem \ref{Main_Theorem}. \hfill\qed

\section{Equivariant Chern classes on $\cQ_d$}\label{classes}

Recall that $\P^1\times\cQ_d$ has a universal exact sequence of sheaves
\[
   0\ \rightarrow\ \cS\ \rightarrow\
   \cO^n_{\P^1\times\cQ_d}\ 
   \rightarrow\  \cT\ \rightarrow 0
\]
with $\cS$ the tautological vector bundle of rank $n-r$.

Since both $\cQ_d$ and $\P^1$ have cell decompositions, we have a
K\"unneth decomposition of Chow rings, 
$A^*(\P^1\times\cQ_d)\simeq A^*(\P^1)\otimes_\Z A^*(\cQ_d)$.
Let $\pi\colon \P^1\times\cQ_d\to \cQ_d$ be the projection. 
For each $1\leq i\leq n{-}r$ we may decompose the Chern class $c_i(\cS)$ 
 \begin{equation}\label{Eq:Chern_Decomposition}
   c_i(\cS)\ =\ \pi^*t_i+h\pi^*u_{i-1}\,,
 \end{equation}
where $t_i,u_i\in A^i(\cQ_d)$ and 
$h$ is the class of a point in $A^1(\P^1)$ pulled back to 
$\P^1\times\cQ_d$.
Note that $u_0=-d$.
Str\o mme~\cite{St87} proved that $A^*(\cQ_d)_\Q$ is generated by the classes 
\[
   \{t_1,\ldots,t_k,u_1,\ldots,u_{k-1}\}\,.
\]

For each $1\leq i\leq n{-}r$, the equivariant Chern
class $c_i^T(\cS)$ localizes at a fixed point $p$ to  
the $i$th elementary symmetric polynomial $e_i$ in the $T$-weights of the
fibre of $\cS$ at $p$.
The fixed points of $\P^1\times\cQ_d$ correspond to $\{0,\infty\}\times \cF$.
For $(\vdelta, \va, \vb)\in\cF$, the bundle $\cS_{(\sdelta, \va, \vb)}$ on
$\P^1$ is a sum of line bundles $\cS_{a_j,b_j}\be_j$ for $\delta_j=1$.
Since $\cS_{a,b}$ has weight $a\bff$ at $0$ and $-b\bff$ at $\infty$, 
the localizations of 
$c_i^T(\cS)$ are  
 \begin{eqnarray*}
  c^T_{i,0,(\sdelta, \va, \vb)}&=&
    e_i(\{\be_j+a_j\bff \mid \delta_j=1\})
  \quad \mbox{at } (0,(\vdelta, \va, \vb))\\
  c^T_{i,\infty,(\sdelta, \va, \vb)}&=&
    e_i(\{\be_j-b_j\bff \mid \delta_j=1\})
  \quad \mbox{at } (\infty,(\vdelta, \va, \vb))\,.
 \end{eqnarray*}

We also have a K\"unneth decomposition in equivariant Chow cohomology:
\[
   A^*_T(\P^1\times\cQ_d)\ \xleftarrow{\ \sim\ }\  A^*_T(\P^1)\otimes_S
   A^*_T(\cQ_d)\,.
\]
To see this, just imitate the argument for ordinary Chow cohomology;
the equivariant Chow cohomology of a variety with an algebraic
cell decomposition will be a free $S$-module, with a module basis
given by the closures of the cells.

An equivariant K\"unneth decomposition of $c^T_i(\cS)$ analogous to 
\eqref{Eq:Chern_Decomposition} requires the choice of 
a lift of the class of a point to $A^1_T(\P^1)$.
Localizing a class $x\in A^1_T(\P^1)$ gives an ordered pair 
$(x_0,x_\infty)\in \Z\bff\oplus\Z\bff$.
Lifts of classes from $A^1(\P^1)$ 
are only well-defined modulo the span of $(\bff,\bff)$.
Three possible choices for lifting the class of a point are:
 \begin{equation}\label{Eq:lifts}
  (i)\  (-\bff,0)\,,\qquad 
  (ii)\  (0,\bff)\,,\quad\mbox{and}\quad
  (iii)\ {\textstyle\frac{1}{2}}(-\bff,\bff)\,.
 \end{equation}
The symmetric lift $(iii)$ requires rational coefficients.
We will use this lift to express our formulas.

Given a lift $h\in A^1_T(\P^1)$ of the class of a point, 
the formula~\eqref{Eq:Chern_Decomposition} defines equivariant lifts of
the classes $t_i,u_{i-1}$.
Str\o mme's result together with Proposition~\ref{P:Chow} implies that 
these classes generate $A^*_T(X)_\Q$ as an $S_\Q$-algebra.

\begin{prop}\label{Prop:kunneth}
 The symmetric choice $(iii)$ of lift $h\in A^1_T(\P^1)$ of the class 
 of a point give the formula for the K\"unneth components $t_i, u_{i-1}$ of the
 equivariant Chern class $c_i^T(\cS)$ of the bundle $\cS$.
 We express this in terms of its localization at the fixed point of $\cQ_d$
 indexed by $p:=(\vdelta, \va, \vb)$.
\[
   t_{i,p}\ =\ \frac{1}{2}(c^T_{i,0,p}+c^T_{i,\infty,p})\qquad
   u_{i-1,p}\ =\ \frac{1}{\bff}(c^T_{i,\infty,p}-c^T_{i,0,p})\,.
\]
\end{prop}

\section{Equivariant Chow ring of the quot scheme
  $\cQ_2(0,2)$}\label{S:ex} 

We use Theorem~\ref{Main_Theorem} to describe the equivariant Chow
ring of $\cQ_2:=\cQ_2(0,2)$.
We first give a basis for $A^*_T(\cQ_2)$ as a module over 
$S=\Z[\be_1,\be_2,\bff]$.

The equivariant Chow ring is the collection of tuples
$(f_p\mid p\in \cQ_2^T)\in S^{ \cQ_2^T}$ which satisfy 
the relations of Theorem~\ref{Main_Theorem}.
If $p$ and $q$ are connected by an edge in the moment multigraph
(see Figure~\ref{F:Q202}) with weight $\chi$,
then $f_p-f_q$ lies in the ideal generated by $\chi$.
These are the standard GKM relations.

There are two
multiedges with four vertices in the multigraph, namely the vertical and horizontal 
lines of symmetry.
They should be seen as flattened quadrangles, since 
they are images of subvarieties isomorphic to $\P^1 \times \P^1$.  
Each 
gives rise to an additional relation, as follows.
Suppose that the quadrangle has four vertices $a,b,c,d$:
\[
  \begin{picture}(150,26)(-50,-1)
   \put(-50,8){\vector(1,0){20}}\put(-45,12){$\chi$}
   \put( 0,10){$a$}\put(42,2){$b$}\put(42,17){$c$}\put(84,10){$d$}
   \put( 9,14){\line(6, 1){30}} \put( 9, 10){\line(6,-1){30}}
   \put(52,19){\line(6,-1){30}} \put(52,  5){\line(6, 1){30}}
  \end{picture}
\]
(Here, the edges are parallel with direction $\chi$.)
Then the tuples $(f_p)$ must satisfy
\[
  f_a - f_b - f_c + f_d \ \in\ \chi^2 S\,.
\]
These are relations of types II(c)$'$ (horizontal) and III$'$ (vertical) of
Theorem~\ref{Main_Theorem}. 

The remaining multiedges with more than two vertices are the left 
and right vertical edges, both with three vertices.  They should
be seen as flattened triangles, since they come from subvarieties
isomorphic to $\mathbb{P}^2$.
The additional relations they induce are described as follows.
Let the three fixed points on the multiedge 
be $a, b, c$, with $b$ between $a$
and $c$.  Then 
\[
  f_a - 2 f_b + f_c\ \in\ \bff^2 S\,.
\]
This is relation of type III$'$  of Theorem~\ref{Main_Theorem}. 

When $k=\C$, we can construct an $S$-module basis for equivariant cohomology
using equivariant Morse theory; as is  well-known, a generic projection of the
moment map to a line will give a Morse function which is perfect for equivariant
cohomology.  
This results in an inductive algorithm to produce a basis, which is nicely
expounded in~\cite{Tym} (see also \cite{GZ01}). 
Pick a direction vector $\bv$ (corresponding to an element of the Lie algebra
$\mathfrak t_\R$) which does not
annihilate the direction vector of any edge.
Orient each edge to have positive pairing with $\bv$, this is the Hasse diagram
of a partial order on the fixed points induced by $\bv$.
Then, using the relations described above, we can inductively construct a 
triangular basis with respect to this ordering.
That is, if $f=f(p)$ corresponds to the fixed point $p$, it vanishes at $q$
($f_q=0$) unless $p < q$, and $f_p$ is the product of weights of
edges pointing down from $p$.
While this algorithm was motivated by Morse theory, 
it makes sense over any field $k$, if $\bv$ is a linear
function on the character group of $T$ which does not annihilate any edge of the
moment graph.

\begin{figure}[htb]
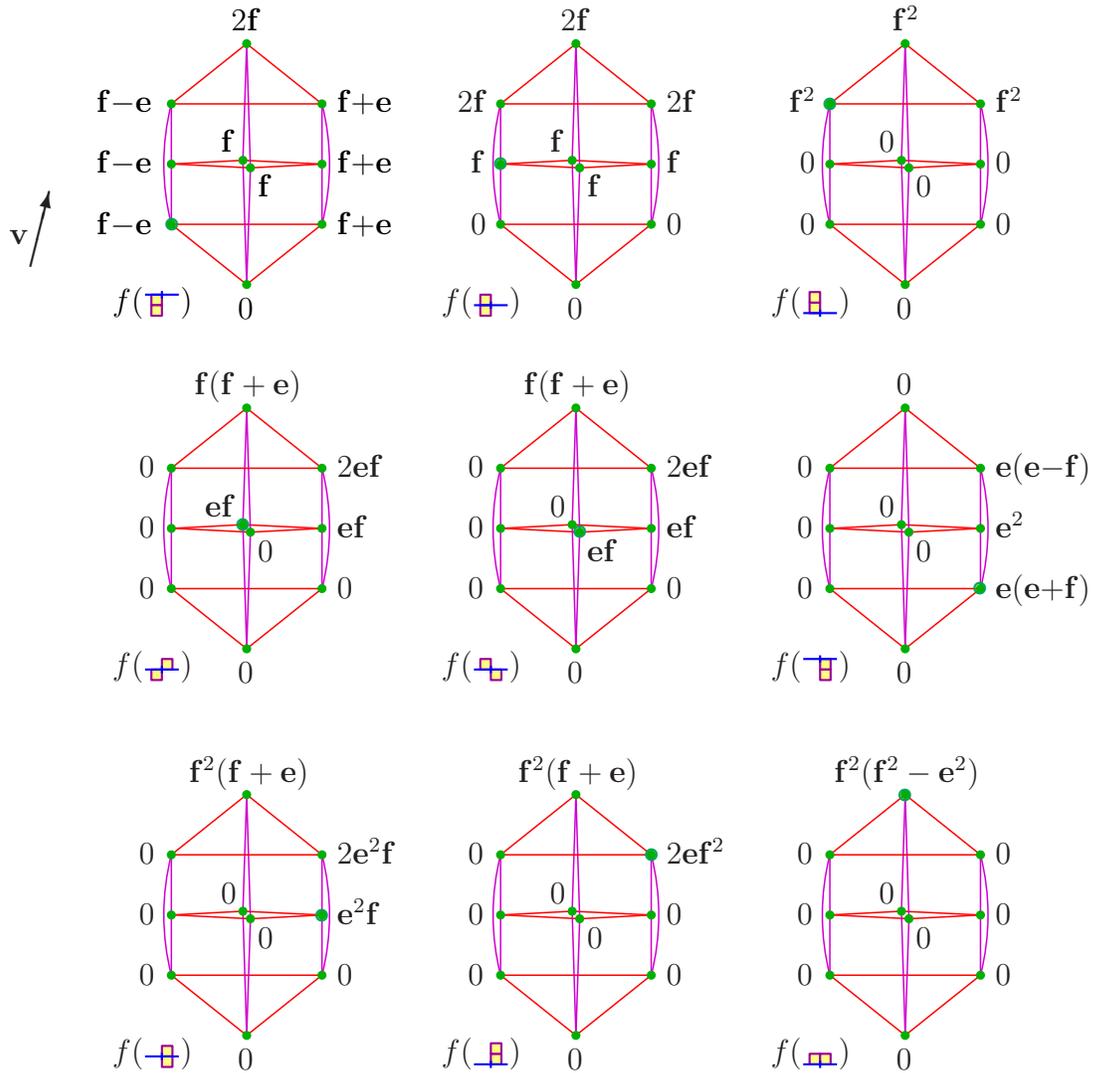

\[
  \begin{picture}(101,130)(-9,-11)
   \put(0,11){\includegraphics[height=95pt]{figures/Q202.skeleton.eps}}
                      \put(26,109){$2\bff$}
     \put(-25,78){$\bff{-}\be$} \put(22,63){$\bff$} \put(66,78){$\bff{+}\be$}
     \put(-25,55){$\bff{-}\be$}                     \put(66,55){$\bff{+}\be$}
     \put(-25,32){$\bff{-}\be$} \put(36,46){$\bff$} \put(66,32){$\bff{+}\be$}
                    \put(28.5,  0){0}
     \put(-19,3){$f(\raisebox{-3pt}{\includegraphics[width=14pt]{figures/11.00.20.eps}})$}
     \thicklines \put(3.4,35.6){\Green{\circle{4}}}
     \put(-50,20){\vector(1,4){7}}
     \put(-58,28){$\bv$}
  \end{picture}
 \qquad
  \begin{picture}(101,130)(-9,-11)
   \put(0,11){\includegraphics[height=95pt]{figures/Q202.skeleton.eps}}
                      \put(26,109){$2\bff$}
     \put(-13,78){$2\bff$}\put(22,63){$\bff$} \put(66,78){$2\bff$}
     \put(-8,55){$\bff$}                      \put(66,55){$\bff$}
     \put(-8,32){0}       \put(36,46){$\bff$} \put(66,32){0}
                    \put(28.5,  0){0}
      \thicklines \put(3.4,58.5){\Green{\circle{4}}}
     \put(-19,3){$f(\raisebox{-3pt}{\includegraphics[width=14pt]{figures/11.10.10.eps}})$}
  \end{picture}
 \qquad
  \begin{picture}(101,130)(-9,-11)
   \put(0,11){\includegraphics[height=95pt]{figures/Q202.skeleton.eps}}
                      \put(27,109){$\bff^2$}
     \put(-12,78){$\bff^2$} \put(22,63){0} \put(66,78){$\bff^2$}
     \put(-8,55){0}                      \put(66,55){0}
     \put(-8,32){0}       \put(36,46){0} \put(66,32){0}
                    \put(28.5, 0){0}
     \thicklines \put(3.4,81.3){\Green{\circle{4}}}
     \put(-19,3){$f(\raisebox{-3pt}{\includegraphics[width=14pt]{figures/11.20.00.eps}})$}
  \end{picture}
\]
\[
  \begin{picture}(101,132)(-9,-11)
   \put(0,11){\includegraphics[height=95pt]{figures/Q202.skeleton.eps}}
                      \put(12,109){$\bff(\bff+\be)$}
     \put(-9,78){0}\put(16,63){$\be\bff$}\put(66,78){$2\be\bff$}
     \put(-9,55){0}                      \put(66,55){$\be\bff$}
     \put(-9,32){0}       \put(36,46){0} \put(66,32){0}
                    \put(28.5,  0){0}
     \thicklines \put(30.3,60){\Green{\circle{4}}}
     \put(-19,3){$f(\raisebox{-3pt}{\includegraphics[width=14pt]{figures/11.01.10.eps}})$}
  \end{picture}
\qquad
  \begin{picture}(101,132)(-9,-11)
   \put(0,11){\includegraphics[height=95pt]{figures/Q202.skeleton.eps}}
                      \put(12,109){$\bff(\bff+\be)$}
     \put(-9,78){0} \put(22,63){0} \put(66,78){$2\be\bff$}
     \put(-9,55){0}                      \put(66,55){$\be\bff$}
     \put(-9,32){0} \put(36,46){$\be\bff$} \put(66,32){0}
                    \put(28.5,  0){0}
      \thicklines \put(33.3,57.2){\Green{\circle{4}}}
    \put(-19,3){$f(\raisebox{-3pt}{\includegraphics[width=14pt]{figures/11.10.01.eps}})$}
  \end{picture}
\qquad
  \begin{picture}(101,132)(-9,-11)
   \put(0,11){\includegraphics[height=95pt]{figures/Q202.skeleton.eps}}
                      \put(28.5,109){0}
     \put(-9,78){0} \put(22,63){0} \put(66,78){$\be(\be{-}\bff)$}
     \put(-9,55){0}                \put(66,55){$\be^2$}
     \put(-9,32){0} \put(36,46){0} \put(66,32){$\be(\be{+}\bff)$}
                    \put(28.5,  0){0}
      \thicklines \put(60.2,35.7){\Green{\circle{4}}}
     \put(-19,3){$f(\raisebox{-3pt}{\includegraphics[width=14pt]{figures/11.00.02.eps}})$}
  \end{picture}
\]
\[
  \begin{picture}(101,132)(-9,0)
   \put(0,11){\includegraphics[height=95pt]{figures/Q202.skeleton.eps}}
                      \put(10,109){$\bff^2(\bff+\be)$}
     \put(-9,78){0} \put(22,63){0} \put(66,78){$2\be^2\bff$}
     \put(-9,55){0}                \put(66,55){$\be^2\bff$}
     \put(-9,32){0} \put(36,46){0} \put(66,32){0}
                    \put(28.5,  0){0}
      \thicklines \put(60.2,58.5){\Green{\circle{4}}}
      \put(-19,3){$f(\raisebox{-3pt}{\includegraphics[width=14pt]{figures/11.01.01.eps}})$}
  \end{picture}
\qquad
  \begin{picture}(101,132)(-9,0)
   \put(0,11){\includegraphics[height=95pt]{figures/Q202.skeleton.eps}}
                      \put(10,109){$\bff^2(\bff+\be)$}
     \put(-9,78){0} \put(22,63){0} \put(66,78){$2\be\bff^2$}
     \put(-9,55){0}                \put(66,55){0}
     \put(-9,32){0} \put(36,46){0} \put(66,32){0}
                    \put(28.5,  0){0}
      \thicklines \put(60.3,81.4){\Green{\circle{4}}}
     \put(-19,3){$f(\raisebox{-3pt}{\includegraphics[width=14pt]{figures/11.02.00.eps}})$}
  \end{picture}
\qquad
  \begin{picture}(101,132)(-9,0)
   \put(0,11){\includegraphics[height=95pt]{figures/Q202.skeleton.eps}}
                      \put(5,109){$\bff^2(\bff^2-\be^2)$}
     \put(-9,78){0} \put(22,63){0} \put(66,78){0}
     \put(-9,55){0}                \put(66,55){0}
     \put(-9,32){0} \put(36,46){0} \put(66,32){0}
                    \put(28.5,  0){0}
     \thicklines \put(31.7,104){\Green{\circle{4}}}
     \put(-19,3){$f(\raisebox{-3pt}{\includegraphics[width=14pt]{figures/11.11.00.eps}})$}
  \end{picture}
\]
\caption{An $S$-module basis for $A^*_T(\mathcal{Q}_2(0,2))$.}
\label{F:module}
\end{figure}

Set $\be:=\be_1-\be_2$ and pick the vector $\bv = \bff +\epsilon \be$, where
$\epsilon>0$ is small.
One basis element is the identity 
$f(\includegraphics[width=14pt]{figures/11.00.11.eps})$,
which localizes to 1 at each fixed point.
We display each of the remaining nine in Figure~\ref{F:module} as a
localization diagram, writing its localizations on a copy of a moment
multigraph. 

\begin{figure}[htb]
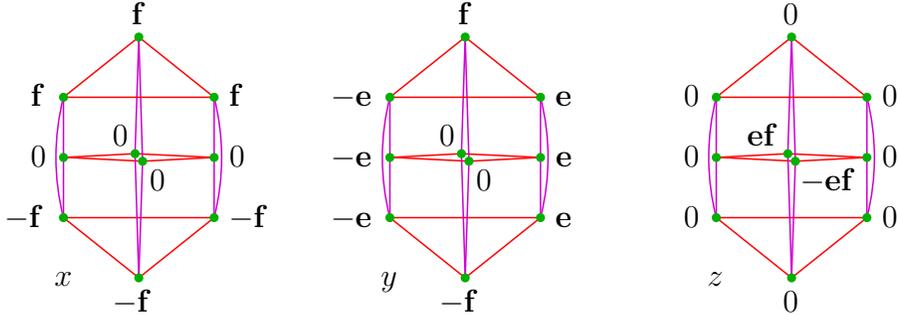

\[
  \begin{picture}(100,130)(-9,0)
   \put(0,11){\includegraphics[height=95pt]{figures/Q202.skeleton.eps}}
                      \put(29,109){$\bff$}
     \put(-9,78){$\bff$} \put(22,63){0} \put(66,78){$\bff$}
     \put(-9,55){0}                \put(66,55){0}
     \put(-19,32){$-\bff$} \put(36,46){0} \put(66,32){$-\bff$}
                    \put(22,  0){$-\bff$}
  \put(0,10){$x$}
  \end{picture}
\qquad
  \begin{picture}(100,130)(-9,0)
   \put(0,11){\includegraphics[height=95pt]{figures/Q202.skeleton.eps}}
                      \put(29,109){$\bff$}
     \put(-19,78){$-\be$} \put(22,63){0} \put(66,78){$\be$}
     \put(-19,55){$-\be$}                \put(66,55){$\be$}
     \put(-19,32){$-\be$} \put(36,46){0} \put(66,32){$\be$}
                    \put(22,  0){$-\bff$}
  \put(0,10){$y$}
  \end{picture}
\qquad
  \begin{picture}(100,130)(-9,0)
   \put(0,11){\includegraphics[height=95pt]{figures/Q202.skeleton.eps}}
                      \put(28.5,109){0}
     \put(-9,78){0} \put(15,63){$\be\bff$} \put(66,78){0}
     \put(-9,55){0}                \put(66,55){0}
     \put(-9,32){0} \put(35,46){$-\be\bff$} \put(66,32){0}
                    \put(28.5,  0){0}
  \put(0,10){$z$}
  \end{picture}
\]
\caption{Generators for the $S_\Q$-algebra $A^*_T(\cQ_2(0,2))_\Q$.}
\label{F:generators}
\end{figure}

Set $x:=f(\raisebox{-3pt}{\includegraphics[width=14pt]{figures/11.10.10.eps}})-
 \bff $, 
$y:=f(\raisebox{-3pt}{\includegraphics[width=14pt]{figures/11.00.20.eps}})-
 \bff $, 
and 
$z:=f(\raisebox{-3pt}{\includegraphics[width=14pt]{figures/11.10.01.eps}})-
    f(\raisebox{-3pt}{\includegraphics[width=14pt]{figures/11.01.10.eps}})$.
Figure~\ref{F:generators} shows their localization diagrams.
We show that they generate $A^*_T(\cQ_2)_\Q$ as an $S_\Q$-algebra
by showing that each basis element $f(p)$ of degree greater than 2
lies in $S_\Q[x,y,z]$. 
Since
 \begin{eqnarray}
    f(\raisebox{-3pt}{\includegraphics[width=14pt]{figures/11.20.00.eps}})
    &=&{\textstyle\frac{1}{2}} x(x{+}\bff)\,,\nonumber\\
    f(\raisebox{-3pt}{\includegraphics[width=14pt]{figures/11.00.02.eps}})
    &=& {\textstyle\frac{1}{2}}(y{+}\be)(y{-}x)\,,\ \
    \mbox{and}\label{Eq:classes}\\
     f(\raisebox{-3pt}{\includegraphics[width=14pt]{figures/11.10.01.eps}})
   + f(\raisebox{-3pt}{\includegraphics[width=14pt]{figures/11.01.10.eps}})\,,
    &=& (y{+}\be)(x{+}\bff)\,,\nonumber\
 \end{eqnarray}
the four degree 2 basis elements lie in $S_\Q[x,y,z]$.
The remaining three basis elements also lie in  $S_\Q[x,y,z]$,
\[
  y f(\raisebox{-3pt}{\includegraphics[width=14pt]{figures/11.01.10.eps}})\ =\ 
    f(\raisebox{-3pt}{\includegraphics[width=14pt]{figures/11.01.01.eps}})\,,
  \quad\quad
  x f(\raisebox{-3pt}{\includegraphics[width=14pt]{figures/11.01.10.eps}})\ =\ 
   f(\raisebox{-3pt}{\includegraphics[width=14pt]{figures/11.02.00.eps}})\,,
  \quad\mbox{and}\quad
   (y-\be) f(\raisebox{-3pt}{\includegraphics[width=14pt]{figures/11.01.01.eps}})\ =\ 
   f(\raisebox{-3pt}{\includegraphics[width=14pt]{figures/11.00.11.eps}})\,.
\]

Inspecting the localization diagrams of $x$, $y$, and $z$, shows that the
following 5 expressions vanish in $A^*_T(\cQ_2)$
 \begin{equation}\label{E:GB}
   xz\,,\quad yz\,,\quad x(x^2-\bff^2)\,,\quad
   (y^2-\be^2)(y-x)\,,\quad \mbox{ and }\quad z^2-(y^2-\be^2)(x^2-\bff^2)\,.
 \end{equation}
In the lexicographic term order where $z>y>x>e>f$, these five polynomials
form a Gr\"obner basis for the ideal ${\mathcal I}$ of 
$\Q[\be_1,\be_2,\bff, \ x, y, z]$ they generate 
with leading terms $ xz$, $yz$, $x^3$, $y^3$, and $z^2$.
There are ten standard monomials 
\[
   1,\ x,\ y,\ x^2,\ xy,\ y^2,\ z,\ x^2y,\ xy^2, \ \mbox{and}\ x^2y^2\,.
\]
Since $A^*_T(\cQ_2)_\Q$ is free over $S_\Q=\Q[\be_1,\be_2,\bff]$ of rank 10, we
conclude that  
\[
    A^*_T(\cQ_2)_\Q \ \simeq\ \Q[\be_1,\be_2,\bff, \ x, y, z]/{\mathcal I}\,.
\] 
Using Proposition~\ref{P:Chow}, we obtain the presentation of the rational Chow
ring 
\[
    A^*(\cQ_2)_\Q \ \simeq\ \Q[x, y, z]/\langle 
    xz,\ yz,\ x^3,\ xy^2-y^3,\ z^2-x^2y^2 \rangle\,.
\] 
By~\eqref{Eq:classes}, the integral Chow ring has a more complicated
presentation.

\begin{figure}[htb]
\[
  \begin{picture}(195,190)(-52,-24)
   \put(-1,11){\includegraphics[height=140pt]{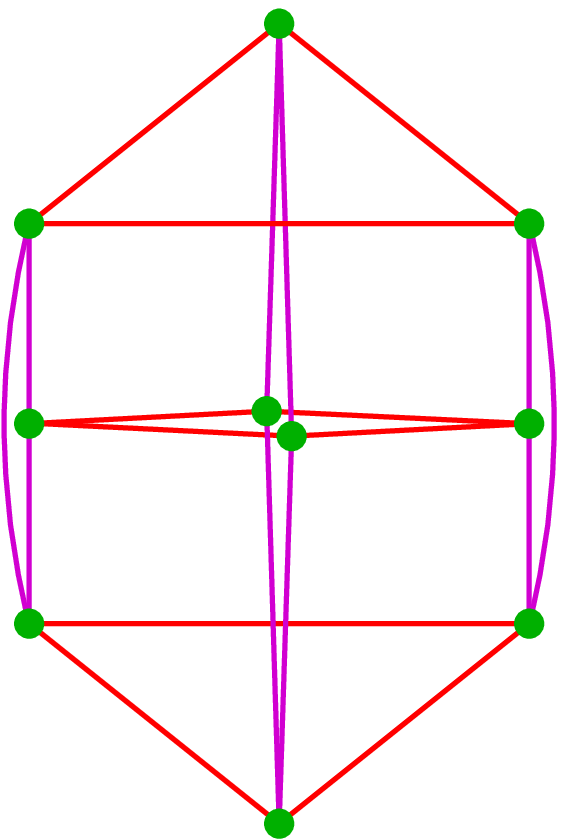}}
    \put(21,155){$\be_1{+}\be_2{+}2\bff$}
    \put(-53,112){$\be_1{+}\be_2{+}2\bff$}
    \put( 93,112){$\be_1{+}\be_2{+}2\bff$} 
    \put(-50, 78){$\be_1{+}\be_2{+}\bff$} 
    \put( 95, 78){$\be_1{+}\be_2{+}\bff$}
    \put(-33, 44){$\be_1{+}\be_2$}
    \put( 93, 44){$\be_1{+}\be_2$}

    \put(10,98){$\be_1{+}\be_2$}
    \put(21,86){${+}\bff$}  
    \put(50,70){$\be_1{+}\be_2$}
    \put(61,58){${+}\bff$} 
    \put(30, 1){$\be_1{+}\be_2$}
    \put(-40,-24){Localization of $c_1^T(\cS)$ at $0\in\P^1$}
  \end{picture}
 \qquad
  \begin{picture}(195,190)(-52,-24)
   \put(-1,11){\includegraphics[height=140pt]{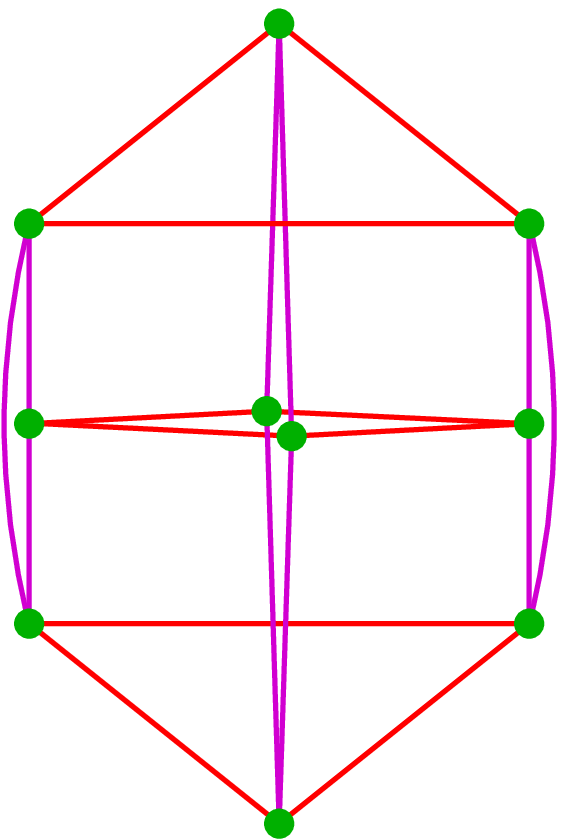}}
    \put( 30,155){$\be_1{+}\be_2$}
    \put(-33,112){$\be_1{+}\be_2$}
    \put( 93,112){$\be_1{+}\be_2$} 
    \put(-50,78){$\be_1{+}\be_2{-}\bff$} 
    \put( 95,78){$\be_1{+}\be_2{-}\bff$}
    \put(-53,44){$\be_1{+}\be_2{-}2\bff$}
    \put( 93,44){$\be_1{+}\be_2{-}2\bff$}

    \put(10,98){$\be_1{+}\be_2$}
    \put(21,86){${-}\bff$}  
    \put(50,70){$\be_1{+}\be_2$}
    \put(61,58){${-}\bff$} 
    \put(21, 1){$\be_1{+}\be_2{-}2\bff$}
    \put(-40,-24){Localization of $c_1^T(\cS)$ at $\infty\in\P^1$}
   \end{picture}
\] 
 \caption{First Chern class of $\cS$.} \label{F:Stromme1}
\end{figure}

We now consider Str\o mme's generators.
Figure~\ref{F:Stromme1} shows the localization diagram of the first Chern class
of $\cS$.
By the formula of Proposition~\ref{Prop:kunneth}, we have
\[
   c^T_1(\cS)\ =\ \be_1+\be_2 + x \,-\, 2h\,,
\]
so that
\[
  t_1\ =\ \be_1+\be_2 + x\qquad\mbox{and}\qquad u_0\ =\ -2\,.
\]

\begin{figure}[htb]
\[
  \begin{picture}(194,190)(-52,-24)
   \put(-1,11){\includegraphics[height=140pt]{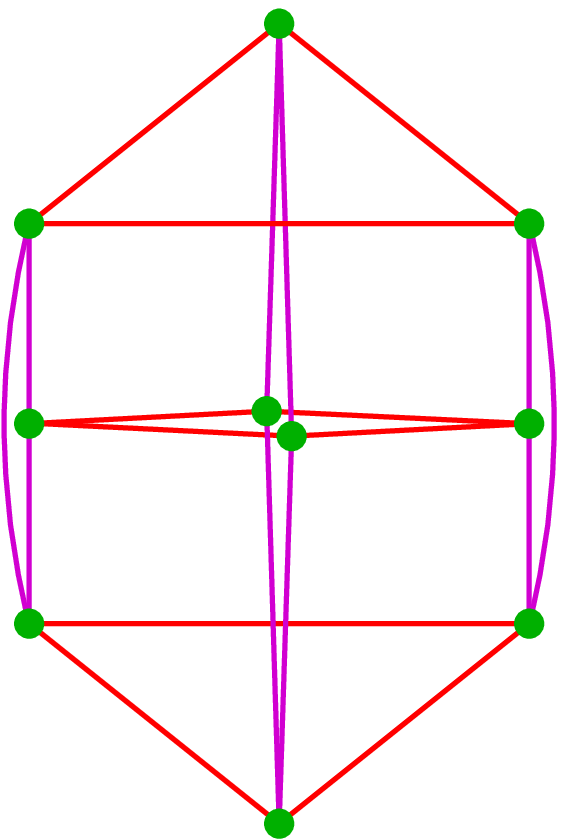}}
    \put(-6,155){$\be_1\be_2{+}(\be_1{+}\be_2)\bff{+}\bff^2$}
    \put(-55,112){$\be_1\be_2{+}2\be_2\bff$}
    \put( 93,112){$\be_1\be_2{+}2\be_1\bff$} 
    \put(-52, 78){$\be_1\be_2{+}\be_2\bff$} 
    \put( 95, 78){$\be_1\be_2{+}\be_1\bff$}
    \put(-25, 44){$\be_1\be_2$}
    \put( 93, 44){$\be_1\be_2$}

    \put(15,100){$\be_1\be_2$}
    \put(15, 88){${+}\be_1\bff$}  
    \put(50, 70){$\be_1\be_2$}
    \put(50, 58){${+}\be_2\bff$} 
    \put(37,  1){$\be_1\be_2$}
    \put(-40,-24){Localization of $c_2^T(\cS)$ at $0\in\P^1$}
  \end{picture}
\qquad
  \begin{picture}(194,190)(-52,-24)
   \put(-1,11){\includegraphics[height=140pt]{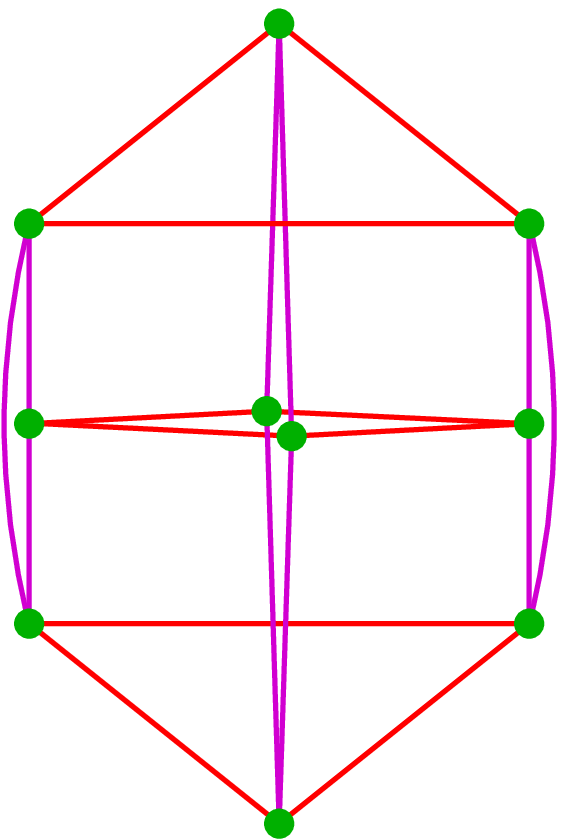}}
    \put( 37,155){$\be_1\be_2$}
    \put(-25,112){$\be_1\be_2$}
    \put( 93,112){$\be_1\be_2$} 
    \put(-52, 78){$\be_1\be_2{-}\be_2\bff$} 
    \put( 95, 78){$\be_1\be_2{-}\be_1\bff$}
    \put(-55, 44){$\be_1\be_2{-}2\be_2\bff$}
    \put( 93, 44){$\be_1\be_2{-}2\be_1\bff$}
    \put(-6,  1){$\be_1\be_2{-}(\be_1{+}\be_2)\bff{+}\bff^2$}

    \put(15,100){$\be_1\be_2$}
    \put(15, 88){${-}\be_2\bff$}  
    \put(50, 70){$\be_1\be_2$}
    \put(50, 58){${-}\be_1\bff$} 
    \put(-40,-24){Localization of $c_2^T(\cS)$ at $\infty\in\P^1$}
  \end{picture}
\]
 \caption{Second Chern class of $\cS$.} \label{F:Stromme2}
\end{figure}

Figure~\ref{F:Stromme2} shows the localization diagram of the second Chern
class of $\cS$.
By the formula of Proposition~\ref{Prop:kunneth}, we have
\[
   c^T_2(\cS)\ =\ \be_1\be_2 + 
    {\textstyle\frac{1}{2}}x(y+ \be_1+\be_2) + {\textstyle\frac{1}{2}}z
     \,-\, h(y+\be_1+\be_2)\,,
\]
so that 
\[
  t_2\ =\ \be_1\be_2 + 
    {\textstyle\frac{1}{2}}x(y+ \be_1+\be_2) + {\textstyle\frac{1}{2}}z
     \qquad\mbox{and}\qquad u_1\ =\ -(y+\be_1+\be_2)\,.
\]

 The corresponding classes in $A^*(\cQ_2)$ are
\[
   t_1\ =\ x\,,\quad
   u_1\ =\ -y\,,\quad\mbox{and}\quad
   t_2\ =\ {\textstyle\frac{1}{2}}(z+xy)\,.
\]
 Note that this shows that the claim in~\cite[Theorem (5.3)]{St87} 
that the classes
 $t_i,u_{i-1}$ generate the integral Chow ring is false.

%
\providecommand{\bysame}{\leavevmode\hbox to3em{\hrulefill}\thinspace}
\providecommand{\MR}{\relax\ifhmode\unskip\space\fi MR }
\providecommand{\MRhref}[2]{%
  \href{http://www.ams.org/mathscinet-getitem?mr=#1}{#2}
}
\providecommand{\href}[2]{#2}

\end{document}